\documentclass[final,1pt,times,twocolumn]{elsarticle}
\usepackage{amssymb}
\usepackage[switch]{lineno}

\usepackage{amsmath}
\usepackage{empheq}%For boxing aligned equations
\usepackage{xcolor}
\usepackage{color, colortbl}
\usepackage[centerlast,small,bf]{caption}
\usepackage{float}
\usepackage{subfig}
\usepackage[T1]{fontenc}
\usepackage[utf8]{inputenc}
\usepackage{appendix}
\usepackage{supertabular}

\newcommand{\myvec}[1]{\mathbf{#1}}

\newcommand{\SetRowColor}[1]{\noalign{\gdef\RowColorName{#1}}\rowcolor{\RowColorName}}
\newcommand{\mymulticolumn}[3]{\multicolumn{#1}{>{\columncolor{\RowColorName}}#2}{#3}}
\definecolor{darkolivegreen}{rgb}{0.33, 0.42, 0.18}
\definecolor{LightBlue}{rgb}{0.9,0.95,1}

\usepackage{tabularx}
\newcolumntype{L}[1]{>{\raggedright\arraybackslash}p{#1}}
\newcolumntype{C}[1]{>{\centering\arraybackslash}p{#1}}
\newcolumntype{R}[1]{>{\raggedleft\arraybackslash}p{#1}}
\newcounter{Eqnno}
\DeclareRobustCommand{\Eqn}[1]{%
   \refstepcounter{Eqnno}%
   \theEqnno\label{#1}}

\journal{Advances in Water Resources}

\begin{document}

\begin{frontmatter}



\title{Multi-rate time stepping schemes for hydro-geomechanical model for subsurface methane hydrate reservoirs}


\author[TUM]{Shubhangi Gupta\corref{cor1}}
\ead{gupta@ma.tum.de}
\cortext[cor1]{We gratefully acknowledge the support for the first author by the German Research Foundation (DFG) through project no. WO 671/11-1.}
\author[TUM,UB]{Barbara Wohlmuth}
\author[IWS]{Rainer Helmig}

\address[TUM]{Chair for Numerical Mathematics, 
              Technical University Munich,
              Boltzmannstra\ss e 3, 85748 Garching bei M\"unchen, Germany    
              }

\address[IWS]{Dept. of Hydromechanics and Modelling of Hydrosystems,
	      University of Stuttgart,
	      Pfaffenwaldring 61, 
	      70569 Stuttgart,
	      Germany 
	      }
	      
\address[UB]{Department of Mathematics, 
              University of Bergen,
              Realfagbygget, All\'egt. 41, Bergen, Norway 
              }	      

\begin{abstract}
We present an extrapolation-based semi-implicit multirate time stepping (MRT) scheme and a compound-fast MRT scheme 
for a naturally partitioned, multi-time-scale hydro-geomechanical hydrate reservoir model.
We evaluate the performance of the two MRT methods compared to an iteratively coupled solution scheme and discuss their advantages and disadvantages.
The performance of the two MRT methods is evaluated in terms of speed-up and accuracy by comparison to an iteratively coupled solution scheme.
We observe that the extrapolation-based semi-implicit method gives a higher speed-up but is strongly dependent on the relative time scales of the latent (slow) and active (fast) components. 
On the other hand, the compound-fast method is more robust and less sensitive to the relative time scales, but gives lower speed up as compared to the semi-implicit method, 
especially when the relative time scales of the active and latent components are comparable.
\end{abstract}


\end{frontmatter}

% \begin{linenumbers}

% \linenumbers

\section{Introduction}\label{sec:introduction}
\paragraph{}
Methane hydrates are crystalline solids formed when water molecules form a cage-like structure and trap a large number of methane molecules within. 
Methane hydrates are thermodynamically stable under conditions of low temperature and high pressure and occur naturally in permafrost regions or sub-seafloor soils. 
If warmed or depressurized, methane hydrates destabilize and dissociate into water and methane gas. 
Natural gas hydrates are considered to be a promising energy resource. 
It is widely believed that the energy content of methane occurring in hydrate form is immense, possibly even exceeding the combined energy content of all other conventional fossil fuels \cite{Kvenvolden1993,Milkov2004}.
Therefore, the development of multiphysics models and numerical codes for coupled hydro-thermo-chemo-geo-mechanical processes are of particular interest  
for evaluating future technologies for gas extraction from methane hydrate reservoirs and for making detailed risk quantification from the inherent geohazards. 

\paragraph{}
A mathematical model describing the hydromechanical processes in a subsurface methane hydrate system has been presented in our earlier work (Gupta et al. (2015)) \cite{Gupta2015}.
The governing PDE's are summarised in Table \ref{table:gverningPDEs} along with some selected closing and constitutive relationships.
This system of PDEs can be decomposed into two sub-classes of models, the \textit{flow and transport model} comprising the mass and energy balance equations 
for the phases occupying the pore spaces in the hydrate formation, i.e., equations \eqref{eqn:MB_CH4}, \eqref{eqn:MB_H2O}, \eqref{eqn:MB_h}, and \eqref{eqn:EB_T}, 
and the \textit{geomechanical model}, comprising the momentum balance equation \eqref{eqn:MomB_sh} for the soil-hydrate composite phase, also referred as the solid-skeleton. 
The model accounts for the effects of the geomechanics on the flow model through the adjustment of the affected reaction and hydraulic properties by scaling the properties with functions of total porosity $\phi$.
In this sense, the soil-phase mass balance (Eqn.\eqref{eqn:MB_s}), which solves for the total porosity, can be seen as mortar between these two sub-models. 
% Since changes in porosity are very small and depend predominantly on soil deformation, 
For simplicity, we eliminate the PDE \eqref{eqn:MB_s} by approximating the total porosity as a function of the volumetric strain 
by assuming that the hydrate-coated soil grains are relatively incompressible compared to the solid skeleton $\epsilon_v$ \cite{Jun-Mo1999}.

\paragraph{}
% The coupling between these two models is dynamic, but weak.
% This observation allows us to decouple the flow and the geomechanical models.
In \cite{Gupta2015}, we have presented an iteratively coupled solution strategy where, the hydrate reservoir model is decomposed into \textit{flow model} and \textit{geomechanical model} as described above,
and is solved iteratively for a given time-step by exchanging shared state variable values between the flow and the geomechaniccal models through a block Gauss-Seidel solution scheme.
At each iteration step, the flow and the geomechanical models solve their corresponding subsystems of equations separately. 
This iteratively coupled solution scheme greatly reduces the computational effort as compared to a monolithic fully implicit scheme.
However, the dynamics of the flow and geomechanical models evolve at different time scales \cite{KimMoridisYangRutqvist2014}. 
We know a priori that the ground deformations manifest at a much slower rate as compared to the flow and transport processes.
Since the refinement of the time-mesh is controlled by the dynamics of active (or the fast) components, 
solving the latent (or the slow) components at this fine time-mesh results in a lot of unnecessary computational work.
The computation can be made cheaper if the slow components are solved on a coarse time-mesh and the fast components on a fine time-mesh.
Such time stepping methods are called \textit{Multi-Rate Time-stepping (MRT) methods}.
The concept of MRT methods was introduced for systems of differential equations (ODEs and DAEs) in such studies as \cite{GearWells1984,GuntherRentrop1993,Rice1960}, 
and some recent results are presented in \cite{BartelGunther2002,SanduConstantinescu2013,EngstlerLubich1997,KvaernoRentrop1999,Savcenco2007}.
MRT methods for hyperbolic conservative laws are developed in \cite{SanduConstantinescu2007,DawsonKirby2001,TangWarnecke2006} and for parabolic equations in \cite{Savcenco2005,Savcenco2008}.
A review of the MRT methods developed over the last two decades can be found in \cite{BookMRTreview}.

\paragraph{}
The application of MRT methods, especially the Implicit-Explicit methods (IMEX), is becoming increasingly popular in the PDE community.
Some of the recent extensions of these methods to application areas of coupled free and porous media flows, air pollution modelling, multi-scale fluid-solid interaction, among others, can be found in 
\cite{RybakHelmigRhode2015,WolkeSchlegel2012,WolkeSchlegel2010,Ehlers2013,ZuijlenBijl2005}.

\paragraph{}
In this article we present two MRT algorithms for our hydro-geomechanical hydrate reservoir model.
The first MRT algorithm is based on a \textit{semi-implicit 'fastest-first'} approach and the second is based on a \textit{compound-fast} approach.
We evaluate the performance of these MRT schemes in comparison to the iteratively coupled solution scheme,
and discuss the advantages and disadvantages of the two MRT methods with respect to the multi-time-scale hydro-geomechanical problems.
To understand the stability of these and related MRT methods in general, the reader is refered to \cite{Verhoeven2007}.

\begin{table*}
  \begin{center}
    \caption{Summary of the mathematical model} \label{table:gverningPDEs}
    \begin{tabular}{L{4cm} R{11cm} L{1cm}}
    \hline\noalign{\smallskip}
    \multicolumn{3}{c}{ \textbf{Governing equations:} }\\[0.5em]
    Mass balance eqn. for each mobile component $\kappa=CH_4,H_2O$  & $\sum\limits_\alpha \partial_t \left( \phi \rho_\alpha S_\alpha \chi_\alpha^\kappa \right)$ 
								      $ + \sum\limits_\alpha \nabla \cdot \left( \phi \rho_\alpha S_\alpha \chi_\alpha^\kappa \myvec{v}_{\alpha,t} \right)$
								      $ = \sum\limits_\alpha \nabla \cdot \left( \phi S_\alpha \myvec{J}_{\alpha}^\kappa \right)$
								      $ + \dot g^\kappa + \sum\limits_\alpha \dot q_{\alpha,m}^\kappa$ 
								    & (\Eqn{eqn:MB_CH4}),(\Eqn{eqn:MB_H2O})\\
    Mass balance eqn. for hydrate phase & $\partial_t \left( \phi \rho_h S_h \right)$  
					  $ + \nabla \cdot \left( \phi \rho_h S_h \myvec{v}_{h,t} \right)$ 
					  $ = \dot g^h $ 
					& (\Eqn{eqn:MB_h})\\
    Mass balance eqn. for soil phase 	& $\partial_t \left[\left(1-\phi\right) \rho_s \right]$  
					  $ + \nabla \cdot \left[\left(1-\phi\right) \rho_s \myvec{v}_{s} \right]$ 
					  $ = 0 $ 
					& (\Eqn{eqn:MB_s})\\
    Energy balance eqn. 		& $\partial_t \left[ \left(1-\phi\right) \rho_s u_s + \sum\limits_\beta \left( \phi \rho_\beta S_\beta u_\beta \right) \right]$ 
					  $ + \sum\limits_\alpha \nabla \cdot \left( \phi \rho_\alpha S_\alpha \chi_\alpha^\kappa \myvec{v}_{\alpha,t} h_\alpha \right)$
					  $ = \nabla \cdot k^c_{eff} \nabla T$
					  $ + \dot Q^h + \sum\limits_\alpha \left( \dot q_{\alpha,m}^\kappa h_\alpha \right) $ 
					& (\Eqn{eqn:EB_T})\\
    momentum balance eqn. for hydrate-soil composite & $ \nabla \cdot \myvec{\tilde\sigma} + \rho_{sh} \myvec{g} = 0 $ 
						     & (\Eqn{eqn:MomB_sh}) \\
    \hline \\[-0.5em]
    \multicolumn{3}{c}{ \textbf{Closing and constitutive relationships:} }\\
    Closure relationships 		& \multicolumn{2}{L{12cm}}{$P_g - P_w = P_c\left(S_{we}\right)$ , \quad $\sum\limits_\beta S_\beta = 1$ , \quad $\forall \kappa$: $\sum\limits_\alpha \chi_\alpha^\kappa = 1$ } \\
    Phase velocities			& \multicolumn{2}{L{12cm}}{	  $\phi S_\beta \myvec{v}_{\beta,t} = \myvec{v}_{\beta} + \phi S_\beta \myvec{v}_{s} $ , 
							    \quad $\myvec{v}_\alpha = \kappa \dfrac{k_{r,\alpha}}{\mu_\alpha} \left( \nabla P_\alpha - \rho_\alpha \myvec{g} \right) $ ,
							    \quad $\myvec{v}_h = 0$ , 
							    \quad $\myvec{v}_s = \partial_t \myvec{u}$} \\
    Diffusive solute flux 		& \multicolumn{2}{L{12cm}}{$\myvec{J}_\alpha^\kappa = - \tau D^\alpha \left( \rho_\alpha \nabla \chi_\alpha^\kappa \right)$ } \\
    Stress-strain relationship		& \multicolumn{2}{L{12cm}}{	  $\myvec{\tilde\sigma} = 2 G_{sh} \myvec{\tilde\epsilon} + \lambda_{sh} \left( tr\ \myvec{\tilde\epsilon} \right) \myvec{\tilde I} + \alpha_{biot} P_{eff}  \myvec{\tilde I}$ , 
							    \quad $\myvec{\tilde\epsilon} = \frac{1}{2}\left( \nabla \myvec{u} + \nabla^T \myvec{u} \right)$}\\
    Reaction kinetics			& \multicolumn{2}{L{12cm}}{ $\dot g^{CH_4} = k^{r} M_{g} A_{rs} \left( P_{eqb}(T) - P_g \right) $ , 
							    \quad $\dot g^{H_2O} = \dfrac{N_h M_{w}}{M_{g}} \dot g^{CH_4} $ , 
							    \quad $\dot g^{h} = \dfrac{M_{h}}{M_{g}} \dot g^{CH_4} $ , 
							    \quad $\dot Q^{h} = \dfrac{\dot g^h}{M_h}\left( B_1 - \dfrac{B_2}{T} \right) $ }\\
    \hline \\[-0.5em]
    \multicolumn{3}{c}{ \textbf{Description:} }\\
    \multicolumn{3}{L{16cm}}{$\alpha = g,w$ \ denotes the mobile gas and liquid phases respectively. The gas phase contains methane gas and water vapour. The liquid phase contains water and dissolved methane.
			     Subscript $s$ denotes the soil phase, $h$ denotes the hydrate phase, and $sh$ denotes the soil-hydrate composite phase.
			     $\beta = g,w,h$ \ denotes all the phases that occupy the pore space.
			     $\kappa = CH_4,H_2O$ \ denotes the molecular components.}\\
    Primary variables &\multicolumn{2}{L{12cm}}{Gas pressure $P_g$, water phase saturation $S_w$, hydrate phase saturation $S_h$, total porosity $\phi$, temprature $T$, and soil displacement $\myvec{u}$. }\\
    Seconday variables &\multicolumn{2}{L{12cm}}{Gas phase saturation $S_g$, water phase pressure $P_w$, mole concentrations $\chi_\alpha^\kappa$, relative phase velocities $\myvec{v}_\beta$, 
						 soil phase velocity $\myvec{v}_s$, total stress $\myvec{\tilde\sigma}$, strain $\myvec{\tilde\epsilon}$. }\\
    Other variables &\multicolumn{2}{L{12cm}}{Methane generation rate $\dot g^{CH_4}$, water generation rate $\dot g^{H_2O}$, hydrate dissociation rate $\dot g^{h}$,
					      heat of dissociation $\dot Q^{h}$, effective pore pressure $P_{eff}$.}\\
    Material properties &\multicolumn{2}{L{12cm}}{Phase density $\rho_\gamma$, capillary pressure $P_c$, intrinsic permeability $\kappa$, relative permeability $k_{r,\alpha}$, phase viscosity $\mu_\alpha$, 
						 binary diffuion cofficient $D^\alpha$, tortuosiy $\tau$, phase entalpy $h_\alpha$, phase internal heat $u_{\gamma}$, effective thermal conductivity $k^c_{eff}$, 
						 lame's parameters $(G_{sh},\lambda_{sh})$, Biot's consant $\alpha_{biot}$,  
						 reaction rate $k^r$, reaction surface area $A_{rs}$, hydrate equilibium pressure $P_{eqb}$, molar mass $M_\beta$, hydration number $N_h$.}\\
    Other material properties &\multicolumn{2}{L{12cm}}{Compressibility $K_\gamma$, Bulk modulus $\gamma$ $B_{\gamma}$, Young's modulus of soil-hydrate composite phase $E_{sh}$,
							Poisson ratio of soil-hydrate composite phase $\nu_{sh}$, thermal conductivities $k^c_{\gamma}$, specific heat capacities $Cp_\alpha,Cv_\gamma$.}\\
    \hline
  \end{tabular}
  \end{center}
\end{table*}

\addtocounter{equation}{\theEqnno}

\section{Multirate time stepping algorithm}\label{sec:MRTmethods}
\paragraph{}
Let the vectors $\myvec{X}_F(t):\mathbb{R}\rightarrow \mathbb{R}^{d_f}$ and  $\myvec{X}_G(t):\mathbb{R}\rightarrow \mathbb{R}^{d_g}$ denote the time-dependent discrete-in-space approximations to the 
primary variables of the \textbf{F}low model (i.e. $P_g,S_w,S_h,T$) and the \textbf{G}eomechanical model (i.e. $\myvec{u}$) respectively (see Table \ref{table:gverningPDEs} for definitions of the variables).
We will refer to $\myvec{X}_F$ as the active components and $\myvec{X}_G$ as the latent components.

\paragraph{}
Further, let $\myvec{F}:\mathbb{R}\times\mathbb{R}^d\times\mathbb{R}^{d_f}\times\mathbb{R}^{d_g}\rightarrow \mathbb{R}^{d_f}$ 
and $\myvec{G}:\mathbb{R}\times\mathbb{R}^d\times\mathbb{R}^{d_f}\times\mathbb{R}^{d_g}\rightarrow \mathbb{R}^{d_g}$
denote the spatial discretization operators for the \textbf{F}low and the \textbf{G}eomechanical models respectively.
Here, $d$ is the dimension of the space domain.
In our numerical scheme, the operator $\myvec{F}$ is obtained by discretizing PDEs (\ref{eqn:MB_CH4}-\ref{eqn:MB_h},\ref{eqn:EB_T}) using the cell-centered finite volume method, 
and the operator $\myvec{G}$ is obtained by discretizing PDE \eqref{eqn:MomB_sh} using the Galerkin finite element method.
However, the MRT methods described in this paper are independent of the methods used for spatial discretization.

\paragraph{}
The spatial discretization of the PDEs (\ref{eqn:MB_CH4}-\ref{eqn:MomB_sh}) governing our hydro-geomechanical model leads to a semi-discrete problem of the following form:

\begin{align}
 &\text{For $t\in[0,T]$, given the initial conditions } \notag\\
 &\quad \myvec{X}_F\left(t=0\right) = \myvec{X}_F^0 \ \text{and, } \myvec{X}_G\left(t=0\right) = \myvec{X}_G^0 \ , \notag\\
 &\text{find solutions for $\myvec{X}_F$ and $\myvec{X}_G$ which satisfy}\notag \\
 &\partial_t \myvec{X}_F = \myvec{F}\left( t, \myvec{x}, \myvec{X}_F, \myvec{X}_G  \right)
 \label{eqn:activeSystem}
 \\
 &\quad \ \ \myvec{0} = \myvec{G}\left( t, \myvec{x}, \myvec{X}_F, \myvec{X}_G  \right)
 \label{eqn:latentSystem}
\end{align}

\paragraph{}
Eqn.\eqref{eqn:activeSystem} is the active ODE (Ordinary Differential Equation) system and Eqn.\eqref{eqn:latentSystem} is the latent AE (Algebraic Equations) system. 
Together, they form a naturally partitioned multi-scale DAE (Differential Algebraic Equations) system.

\begin{figure}
 \centering
  \includegraphics[scale=0.65]{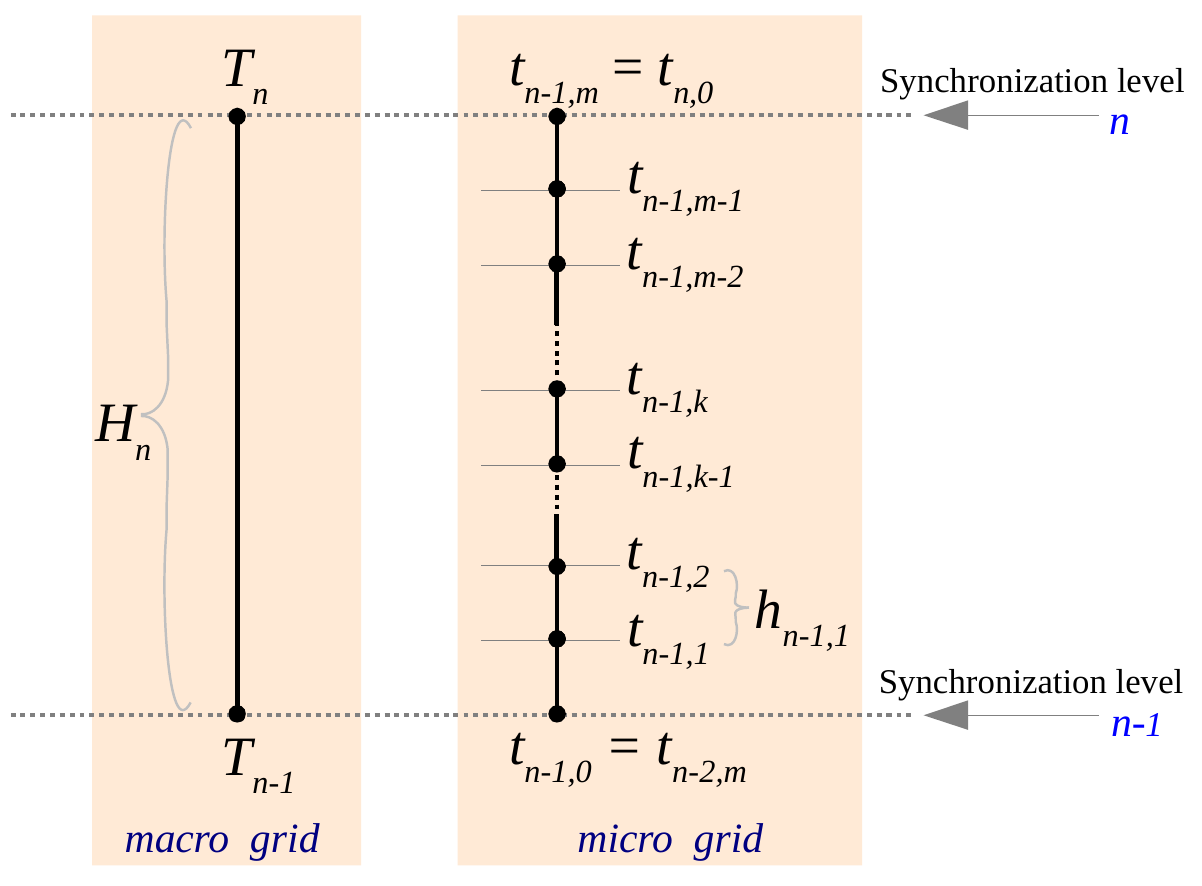}
  \caption{Time-mesh for active and latent components}
  \label{fig:timeMesh}
\end{figure}

\paragraph{}
Each part of the partitioned DAE system is marched in time on an independent time-mesh which depends on it's own \textit{activity}.
Here, activity of a component refers to the time-scale at which the dynamics of the governing equation for that component evolves.
We assume that the activity of the components does not vary in space, i.e., that the components evolve on the same time-scale throughout the spatial domain.
For the latent system, we define a coarse time-mesh $\{ T_n, 0 \leq n \leq N \}$ with time step sizes $\{H_n = T_n - T_{n-1}, 0 < n \leq N \}$. 
We will refer to this as the macro-grid, and the time step from $T_{n-1}$ to $T_n$ as the macro-step.
For the active system, we define a refined time-mesh $\{ t_{n,k} \ , 0 \leq n < N, 0 \leq k \leq m  \}$ with time step sizes $\{ h_{n,k} = t_{n,k} - t_{n,k-1}, 0 \leq n < N, 0 < k \leq m \}$
and multirate factor $m$.
We will refer to this as the micro-grid, and the time step from $t_{n-1,k-1}$ to $t_{n-1,k}$ for each $k=[1,...,m]$ as the micro-steps.
The two time-meshes are synchronized, which implies that for all n, $T_n = t_{n,0} = t_{n-1,m}$ (See Fig. \ref{fig:timeMesh}).

\paragraph{}
All MRT methods have the basic property that the time integration can proceed from synchronization level $n$ to $n+1$ only when all the components, slow and fast, have made their resepective macro and micro steps
and have synchronized at the level $n$. 

\paragraph{}
For marching the DAE system (\ref{eqn:activeSystem},\ref{eqn:latentSystem}) forward in time from $T_{n-1}$ to $T_n$, 
the active ODE \eqref{eqn:activeSystem} is integrated on the micro grid $t_{n-1,k}$ using the implicit Euler method for each micro step,
while the latent AE \eqref{eqn:latentSystem} is evaluated directly at the macro grid point $T_n$ using the solution of the active ODE at $t_{n-1,m}$.
The two MRT algorithms that we will discuss differ in how the latent components are approximated on the micro grid for solving the active ODE.

\paragraph{}
In the semi-implicit MRT method (Algorithm 1), we first make the $m$ micro steps for the active components from $t_{n-1,0}$ to $t_{n-1,m}$.
The values of the latent components needed on the micro grid, i.e. ${\myvec{X}_G}_{n-1,k}$, for making the micro steps are approximated by means of extrapolation.
In our scheme, we construct a polynomial function of order $p$ for extrapolation using the values of $\myvec{X}_G$ evaluated at $p+1$ previous macro grid points, 
i.e., $T_{n-1},...,T_{n-(p+1)}$.
We then make the final macro step to evaluate the latent component at $T_n$.

\paragraph{}
For $m=1$, this method essentially becomes a decoupled sequential solution scheme, which by itself is faster than the iteratively coupled solution scheme.
For $m\geq 1$, all systems are solved only once on their respective time-meshes, 
and the coefficients of the extrapolation function also need evaluation only once per macro step, 
thus requiring very little computational effort.

\begin{figure}
  \begin{center}\textcolor{gray}{\rule{8.cm}{1pt}}\end{center}
  \vspace{-35pt}
  \paragraph{ALGORITHM 1: Semi-implicit MRT method}
  \begin{itemize}
  \item[] \textbf{STEP 1: Extrapolation macro step} \newline Extrapolate ${\myvec{X}_G}_{n-1,k}$ at each $k = [1...m]$ on the fine time-mesh using the $p+1$ old step values of ${\myvec{X}_G}$ 
      computed at $T_{n-1},...,T_{n-(p+1)}$:
      \begin{align}\label{eqn:alg1-extrapolation}
      &{\myvec{\tilde{X}{}}_{G}}_{n-1,k} = {\myvec{X}_{G}}_{n-1} + \sum \limits_{j=1}^{p} A_j \left( t_{n-1,k} - T_{n-1} \right)^{j} 
%       \\
%       &\text{where, }
%       \notag
%       \\
%       &A_2 = \frac{1}{H_{n-2}}\left(\frac{ {\myvec{X}_G}_{n-1} - {\myvec{X}_G}_{n-2} }{ H_{n-1} } - \frac{ {\myvec{X}_G}_{n-1} - {\myvec{X}_G}_{n-3} }{ H_{n-1} + H_{n-2} } \right)\notag
%       \\
%       &A_1 = \frac{ {\myvec{X}_G}_{n-1} - {\myvec{X}_G}_{n-2} }{ H_{n-1} } + A_2 \ H_{n-1} \notag
      %  \left( {\myvec{X}_G}_{n-1}, {\myvec{X}_G}_{n-2}, {\myvec{X}_G}_{n-3} \right)
      \end{align}
    \item[] \textbf{STEP 2: Micro-steps} \newline Solve for ${\myvec{X}_F}_{n-1,k}$ at each $k = [1,...,m]$ on the fine time-mesh using the implicit Euler method:
      \begin{align}\label{alg:1-2}
      &{\myvec{X}_F}_{n-1,k} = {\myvec{X}_F}_{n-1,k-1} \notag \\
			  &+ h_{n-1,k} \ \myvec{F}\left( t_{n-1,k} \ , \ \myvec{x} \ , \ {\myvec{X}_F}_{n-1,k} \ , \ {\myvec{\tilde{X}{}}_G}_{n-1,k} \right) 
      \end{align}
    \item[] \textbf{STEP 3: Macro-step} \newline Solve for ${\myvec{X}_G}_{n}$:
      \begin{align}
      \myvec{G} \left( T_{n} \ , \ \myvec{x} \ , \ {\myvec{X}_F}_{n-1,m}\ , \ {\myvec{X}_G}_{n} \right) = 0.
      \end{align}
  \end{itemize}
  \vspace{-25pt}
  \begin{center}\textcolor{gray}{\rule{8.cm}{1pt}}\end{center}
\end{figure}

\paragraph{}
In the compound-fast MRT method (Algorithm 2), we first make a predictor macro step to get an approximate value of the latent component at $T_n$.
In this step, we integrate the active ODE on the macro grid from $T_{n-1}$ to $T_n$ with a relaxed stopping criteria for the Newton slover. 
This gives a rough approximation of $\myvec{X}_F$ at $T_n$ (denoted by ${\myvec{\tilde{X}{}}_F}_{n}$), 
which is then used to solve the latent AE to get an approximate value of $\myvec{X}_G$ at $T_n$ (denoted by ${\myvec{\tilde{X}{}}_G}_{n}$).
We then make the micro steps to integrate the active ODE from $t_{n-1,0}$ to $t_{n-1,m}$.
The values of the latent components needed on the micro grid, i.e. ${\myvec{X}_G}_{n-1,k}$, for making the micro steps are approximated by means of linear interpolation (refer Eqn.\eqref{eqn:alg2-interpolation}).
In the final step, called the corrector macro step, we solve the latent AE once more at the macro grid point $T_n$ to correct (improve) the solution from the predictor step.

\paragraph{}
If the ODE becomes unsolvable on the macro grid and the predictor step fails, then, in our simulator we reduce the value of $H_n$ by half and attempt the predictor step once again.

\begin{figure}
  \begin{center}\textcolor{gray}{\rule{8.cm}{1pt}}\end{center}
  \vspace{-35pt}
  \paragraph{ALGORITHM 2: Compound-fast MRT method}
  \begin{itemize}
  \item[] \textbf{STEP 1: Predictor macro step} (or compound step) 
      \newline Integrate active ODE \eqref{eqn:activeSystem} with large step size $H_n$. Relax the stopping criteria for the Newton solver to get a rough approximation at $T_n$, ${\myvec{\tilde{X}{}}_F}_{n}$: 
      \begin{align}\label{alg:2-1}
      {\myvec{\tilde{X}{}}_F}_{n} &= {\myvec{X}_F}_{n-1} + H_{n} \cdot \myvec{F}\left( T_{n} \ , \ \myvec{x} \ , \ {\myvec{\tilde{X}{}}_F}_{n} \ , \ {\myvec{X}_G}_{n-1} \right) 
      \end{align}
      Use ${\myvec{\tilde{X}{}}_F}_{n}$ to predict ${\myvec{\tilde{X}{}}_G}_{n}$:
      \begin{align}
      \myvec{G} \left( T_{n} \ , \ \myvec{x} \ , \ {\myvec{\tilde{X}{}}_F}_{n}\ , \ {\myvec{\tilde{X}{}}_G}_{n} \right) = 0.
      \end{align}

    \item[] \textbf{STEP 2: Micro-steps} 
      \newline Solve for ${\myvec{X}_F}_{n-1,k}$ at each $k = [1,...,m]$ on the fine time-mesh using implicit Euler method:
      \begin{align}\label{alg:2-2}
      &{\myvec{X}_F}_{n-1,k} = {\myvec{X}_F}_{n-1,k-1} \notag \\
			  &+ h_{n-1,k} \ \myvec{F}\left( t_{n-1,k} \ , \ \myvec{x} \ , \ {\myvec{X}_F}_{n-1,k} \ , \ {\myvec{\tilde{X}{}}_G}_{n-1,k} \right) 
      \end{align}
      where, ${\myvec{\tilde{X}{}}_G}_{n-1,k}$ are the linearly interpolated values of $\myvec{X}_G$ at $t_{n-1,k}$:
      \begin{align}\label{eqn:alg2-interpolation}
      {\myvec{\tilde{X}{}}_G}_{n-1,k} = {\myvec{X}_G}_{n-1} + \left( {\myvec{\tilde{X}{}}_G}_{n} - {\myvec{X}_G}_{n-1} \right) \sum\limits _{i=1}^k \frac{h_{n-1,k}}{H_n}
      \end{align}

    \item[] \textbf{STEP 3: Corrector macro step} 
      \newline Solve for ${\myvec{X}_G}_{n}$:
      \begin{align}
      \myvec{G} \left( T_{n} \ , \ \myvec{x} \ , \ {\myvec{X}_F}_{n-1,m}\ , \ {\myvec{\tilde{X}{}}_{G}}_{n} + \Delta {\myvec{X}_{G}}_{n} \right) = 0.
      \end{align}
  \end{itemize}
  \vspace{-30pt}
  \begin{center}\textcolor{gray}{\rule{8.cm}{1pt}}\end{center}
\end{figure}

% \paragraph{Adaptive time-step control: }
% For problems where an adaptive time step control is desired to increase the speed up, especially where the activity of the latent components relative to that of the active components is expected to change over time, 
% we can control the refinement of the macro grid in the compound-fast MRT method.
% We use the number of newton steps $i_{p,n}$ required by the active system in the predictor step from level $n-1$ to level $n$ to heuristically adjust the next macro-step size $H_{n+1}$ (Refer Algorithm 3).
% \begin{figure}
%  \begin{center}\textcolor{gray}{\rule{8.4cm}{1pt}}\end{center}
%     \vspace{-15pt}
%   \paragraph{ALGORITHM 3: Adaptive macro-step control for compound-fast MRT method }
%    \begin{itemize}
%     \item[] If $\ i_{p,n}\leq 2\ $, then $\ H_{n+1}=a \cdot H_{n}\ $, where $\ a>1\ $.\newline
%     We ususally choose $\ 1< a \leq 2\ $. \vspace{0.2cm}
%     \item[] Else if $\ i_{p,n}\geq 6\ $, then $\ H_{n+1}=b \cdot H_{n}\ $, where $\ 0<b<1\ $.\newline
%     We ususally choose $\ 0.75\leq b \leq 0.9\ $.  \vspace{0.2cm}
%     \item[] Else, $\ H_{n+1}=H_{n}\ $.
%    \end{itemize}
%     \vspace{-20pt}
%   \begin{center}\textcolor{gray}{\rule{8.4cm}{1pt}}\end{center}
% \end{figure}

\paragraph{}
Both the MRT algorithms are implemented in the C++ based DUNE PDELab framework \cite{DUNEBastian2007} as an extension to our hydro-geomechanical hydrate reservoir simulator.
This code is capable of solving problems in 1D, 2D and 3D domains. The MRT methods discussed in this section are however, independent of the dimension of the space domain.
% For sake of simplicity, in this paper we will present numerical results only for a 1D test problem.

\section{Numerical examples}\label{sec:numericalExamples}
\paragraph{}
We now present two model problems to test the performance of the MRT algorithms presented in Section \ref{sec:MRTmethods}.
The first problem considers $1D$ consolidation in a depressurized methane hydrate sample. 
The second problem considers a relatively complex $3D$ example where we simulate the hydro-geomechanical processes in a subsurface hydrate reservoir 
which is destabilized by depressurization using a vertically placed low pressure gas well.

\subsection{Test 1}
\paragraph{}
We consider a 1D consolidation in a depressurized methane hydrate sample as our test problem.

\subsubsection{Problem setting}\label{sec:problemSetting}
\begin{figure}
 \centering
 \includegraphics[scale=0.4]{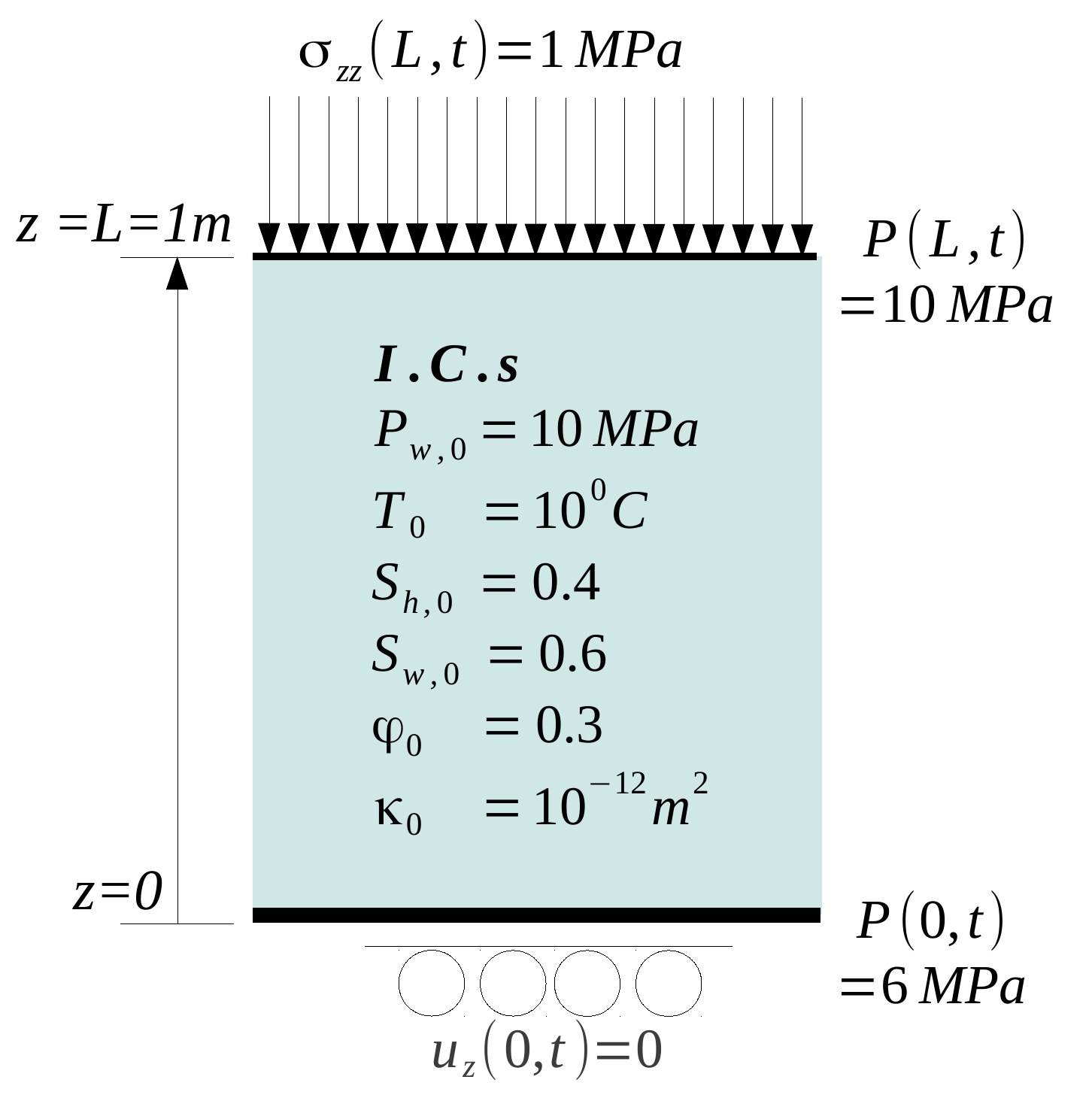}
 \caption{Schematic of the $1D$ test problem}
 \label{fig:problemSetting}
\end{figure}

\paragraph{}
The problem set-up consists of a confined soil sample of height $1$m. 
The sample has a uniformly distributed hydrate saturation of $S_h=0.4$ and is fully saturated with water. 
The porosity and the permeability of the hydrate free soil are $\phi=0.3$ and $\kappa=10^{-12}$ $m^2$ respectively.
A constant vertical stress of $1$ MPa acts at the top boundary while the lower boundary is held fixed. 
At the upper boundary the pressure is kept constant at the initial value of $10$ MPa, while at the lower boundary a low pressure of $6$ MPa is maintained at all times. 
The problem setting imitates gas production in a hydrate reservoir through depressurization which is achieved by pumping out the free gas occuring below the gas-hydrate layer.
The schematic for this problem is shown in Fig. \ref{fig:problemSetting}.
The selected material properties and model parameters are listed in Table \ref{table:properties}.

\begin{table}
\caption{Material properties and model parameters for the $1D$ test problem} \label{table:properties}
\begin{tabular}{L{0.5cm} C{5.25cm} L{1cm}}
\hline
      \SetRowColor{LightBlue}
      \mymulticolumn{3}{c}{\textbf{Thermal conductivities}} \\ \\[-0.9em]
      $k^c_g$				& $- 0.886\times 10^{-2}\ \ \ 
					  + 0.242\times 10^{-3}\ T\ 
					  - 0.699\times 10^{-6}\ T^2 
					  + 0.122 \times 10^{-8}\ T^3$	& $\frac{W}{m\cdot K}$  \\
      $k^c_w$				& $0.3834\ ln(T) - 1.581$	& $\frac{W}{m\cdot K}$	\\
      $k^c_h$				& $2.1$				& $\frac{W}{m\cdot K}$	\\
      $k^c_s$				& $1.9$				& $\frac{W}{m\cdot K}$	\\ \\[-0.8em]
% \\
%       
      \SetRowColor{LightBlue}
      \mymulticolumn{3}{c}{\textbf{Specific heat capacities}} \\ \\[-0.9em]
      $Cp_w$				& $4186$		& $\frac{J}{kg\cdot K}$ \\
      $Cv_w$				& $Cp_w+R_{H_2O}$	& $\frac{J}{kg\cdot K}$ \\
      $Cv_h$				& $2700$		& $\frac{J}{kg\cdot K}$ \\
      $Cv_s$				& $800$			& $\frac{J}{kg\cdot K}$ \\ \\[-0.8em]
% \\
%       
      \SetRowColor{LightBlue}
      \mymulticolumn{3}{c}{\textbf{Dynamic viscosities}} \\ \\[-0.9em]
      $\mu_g$	& $10.4 \ e{-6} \left(\dfrac{273.15 + 162 }{T+162}\right) \left(\dfrac{T}{273.15}\right)^{1.5}$	& $Pa\cdot s$	\\
      $\mu_w$	& $0.001792 \ \exp\left[ - 1.94 - 4.80 \left(\dfrac{273.15}{T}\right) \right.$			& $Pa\cdot s$	\\
		& $\left. + 6.74 \left(\dfrac{273.15}{T}\right)^2\right]$ 					& 		\\ \\[-0.8em]
% \\
% 
      \SetRowColor{LightBlue}
      \mymulticolumn{3}{c}{\textbf{Densities}} \\ \\[-0.9em]
      $\rho_g$	& $\dfrac{P_g}{z R_g T}$	& $\frac{kg}{m^{3}}$ 	\\
      $\rho_w$	& $1000$			& $\frac{kg}{m^{3}}$ 	\\
      $\rho_h$	& $900$ 			& $\frac{kg}{m^{3}}$ 	\\
      $\rho_s$	& $2100$ 			& $\frac{kg}{m^{3}}$ 	\\ \\[-0.8em]
% \\
% 
      \SetRowColor{LightBlue}
      \mymulticolumn{3}{c}{\textbf{Hydraulic properties}} \\ \\[-0.9em]
      \multicolumn{3}{l}{Brooks-Corey parameters}\\
      $\lambda_{BC}$		& $1.2$		&	\\
      $P_{entry}$		& $50$		& $kPa$	\\ \\[-0.8em]
% \\
%      
      \SetRowColor{LightBlue}
      \mymulticolumn{3}{c}{\textbf{Hydrate kinetics}} \\ \\[-0.9em]
      $k^{r}$		&$3.6 \times 10^{4} \ \exp\left( - \dfrac{9752.73}{T}\right)$	& $\frac{\text{mol}}{m^{2}\cdot Pa \cdot s} $\\
      $N_{h}$		& $5.75$							& \\
      $P_{eqb}$	 	& $ \exp\left(14.17 - \dfrac{1886.79}{T}\right)$		& $MPa$ \\
      $\dot Q_h$ 	& $A_1=56599$, $A_2=16.744$					& $\frac{W}{m^{3}}$ \\ \\[-0.8em]
% \\
% 
      \SetRowColor{LightBlue}
      \mymulticolumn{3}{c}{\textbf{Poroelasticity parameters}} \\ \\[-0.9em]
      $\alpha_{biot}$	& $0.8$		& \\
      $\nu_{sh}$	& $0.15$	& \\
      $E_{sh}$		&$160 + 250 \ S_h $		& $MPa$\\\hline
\end{tabular}
\end{table}

\subsubsection{Numerical simulation}\label{sec:numSimulation}
\paragraph{}
We discretize the space domain into $200$ cells along the Z-axis, and run the simulation up to $T=18000$ sec.
For simplicity, we assume uniform, non-adaptive micro and macro time-grids. 
The time step size for the micro grid is fixed at $h_{n,k}=h=60$ sec for all $n$ and $k$. 
This time-step satisfies the CFL condition for the flow model to ensure stability of the active system.
We chose different values of the multirate factor, as, 

$$ m=[1,2,5,10,20,30], $$ 

so that the macro grid is $m$ times coarser than the micro grid, i.e. $H_n = H = mh$ for all $n$.

% \begin{figure*}
%  \centering
%  \subfloat[$P_g$ profile over $z$ at $t=3600$ sec.]{
% 	      \includegraphics[scale=0.265]{FC-Pg_overZ_t=3600s.jpg} \label{fig:FC-Pg-over-z-wrt-m}
% 	    }
%   \hspace{-0.25cm}
%   \subfloat[$S_h$ profile over $z$ at $t=3600$ sec.]{
% 	      \includegraphics[scale=0.265]{FC-Sh_overZ_t=3600s.jpg} \label{fig:FC-Sh-over-z-wrt-m}
% 	    }
%   \hspace{-0.25cm}
%   \subfloat[$u_z$ profile over $z$ at $t=3600$ sec.]{
% 	      \includegraphics[scale=0.265]{FC-uz_overZ_t=3600s.jpg} \label{fig:FC-uz-over-z-wrt-m}
% 	    }
%   \caption{$P_g$, $S_h$ and $u_z$ profiles for the $1D$ test problem computed using fully-coupled (FC) scheme.}
%   \label{fig:FC-PgUz-over-z-wrt-m}
% \end{figure*}

\paragraph{}
We evaluate the performance of the two MRT methods in terms of 1) Speed up, and 2) Relative error.
Speed up is calulated as,

\begin{align*}
 \text{speed-up} = \frac{\text{CPU-time for $m$-step MRT method}}{\text{CPU-time for iteratively coupled scheme}} 
\end{align*}

The time step size of $h=60$ sec is used for the iteratively coupled scheme. 
\newline
The relative error is calculated as,

\begin{align*}
 \text{relative error} = \frac{\text{L2-error by $m$-step MRT method}}{\text{L2-error by iteratively coupled scheme}}
\end{align*}

where, to compute the error in the solution from the iteratively coupled scheme and the MRT schemes, the solution from a fully coupled fully implicit scheme is used as the reference solution.
The time step size used for the fully coupled fully implicit scheme is also $h=60$ sec.

\paragraph{}
For the solution of the active system on the micro grid, a minimum error reduction of $10^{-8}$ is prescribed for the Newton solver.
For the predictor step of the Compound-fast MRT method, the minimum error reduction for the Newton solver is relaxed to $10^{-3}$.

\begin{figure}
 \centering
 \subfloat[$P_g$ profile over $z$.]{
	    \includegraphics[scale=0.365]{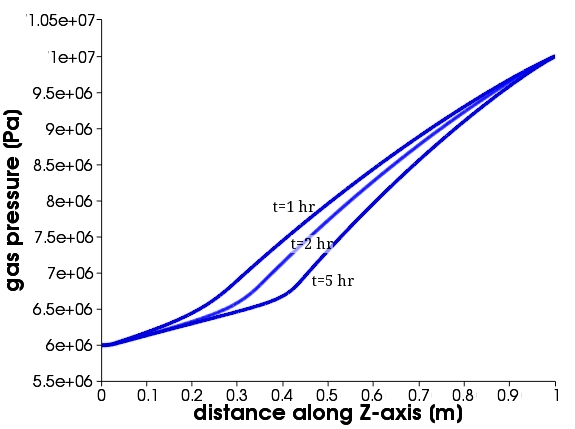} \label{fig:FC-Pg-over-z-wrt-m}
	    }
%   \vspace{0.25cm}
  \vfill
  \subfloat[$u_z$ profile over $z$.]{
	      \includegraphics[scale=0.365]{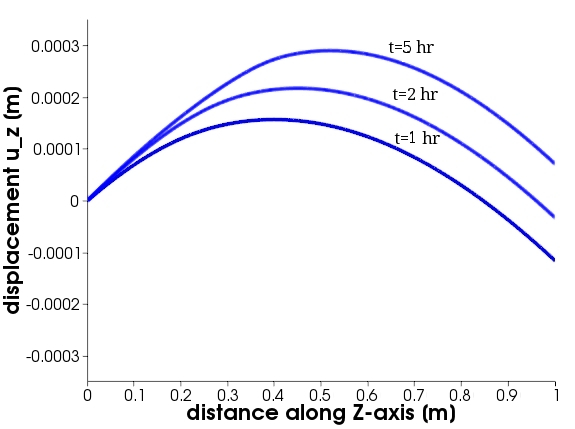} \label{fig:FC-uz-over-z-wrt-m}
	    }
  \caption{$P_g$, $S_h$ and $u_z$ profiles for the $1D$ test problem computed using fully-coupled (FC) scheme.}
  \label{fig:FC-PgUz-over-z-wrt-m}
\end{figure}

\paragraph{}
For reference, the solution of the problem computed using the fully coupled fully implicit scheme is shown in Fig. \ref{fig:FC-PgUz-over-z-wrt-m}.

\subsubsection{Results - Semi-implicit MRT method}
% \paragraph{Semi-implicit MRT method: }
% \paragraph{}
% Fig. \ref{fig:ALG1-Pg-over-z-wrt-m} and Fig. \ref{fig:ALG1-uz-over-z-wrt-m} show snap-shots of gas-phase pressure $P_g$ profiles and the vertical displacement $u_z$ profiles respectively at $t=3600$ sec
% for the fully coupled, iteratively coupled and the m-step semi-implicit MRT schemes. 
% In Fig. \ref{fig:ALG1-speedup-L2errPg-L2errU-vs-m} the relative errors in $P_g$ and $u_z$, and the speed up are plotted over $m$.
% We can see that the errors grow exponentially with increasing $m$, while the speed up shows a steep increase in the beginning and eventually leads to a plateau as $m$ increases. 
% This implies that we cannot choose any arbitrary value of $m$. 
% $m$ should be small enough so that the extrapolation error does not dominate, but it should be large enough so that we can take advantage of the speed up optimally, thus making the choice of $m$ an important consideration. 

% \begin{figure}
%  \centering
%   \includegraphics[scale=0.5]{ALG1-speedup_errorPg_errorU_vs_m.pdf}
%   \caption{a) Speed-up, \ b) relative error in $P_g$, and \ c) relative error in $u_z$ \ over $m$ at $t=18000$ sec for the semi-implicit fast-first MRT method.}
%   \label{fig:ALG1-speedup-L2errPg-L2errU-vs-m}
% \end{figure}

\subsubsection*{Relative error: }
In Fig. \ref{fig:ALG1-Pg-over-m-P0,P1,P2,P3ex}, the relative error in $P_g$ is plotted over $m$ for the semi-implicit MRT method using polynomials of order $p=0,1,2,3$ in the extrapolation steps.
The scheme with polynomial of order $p=0$ is the most stable, but gives only $\approx 80\%$ accuracy as compared to the iteratively coupled solution scheme,
while the polynomial of order $p=2$ gives an accuracy of $\approx 99\%$ compared to the iteratively coupled scheme for $m\leq 5$, but becomes increasingly unstable for higher $m$.
The polynomial of order $p=3$ does not give any significant advantage in terms of accuracy but makes the scheme highly unstable.
% Thus, we can improve the accuracy of the solution by increasing the order of the extrapolation polynomial we can improve the accuracy of the solution at the expense of stability with respect to $m$.

\begin{figure}
 \centering
  \includegraphics[scale=0.5]{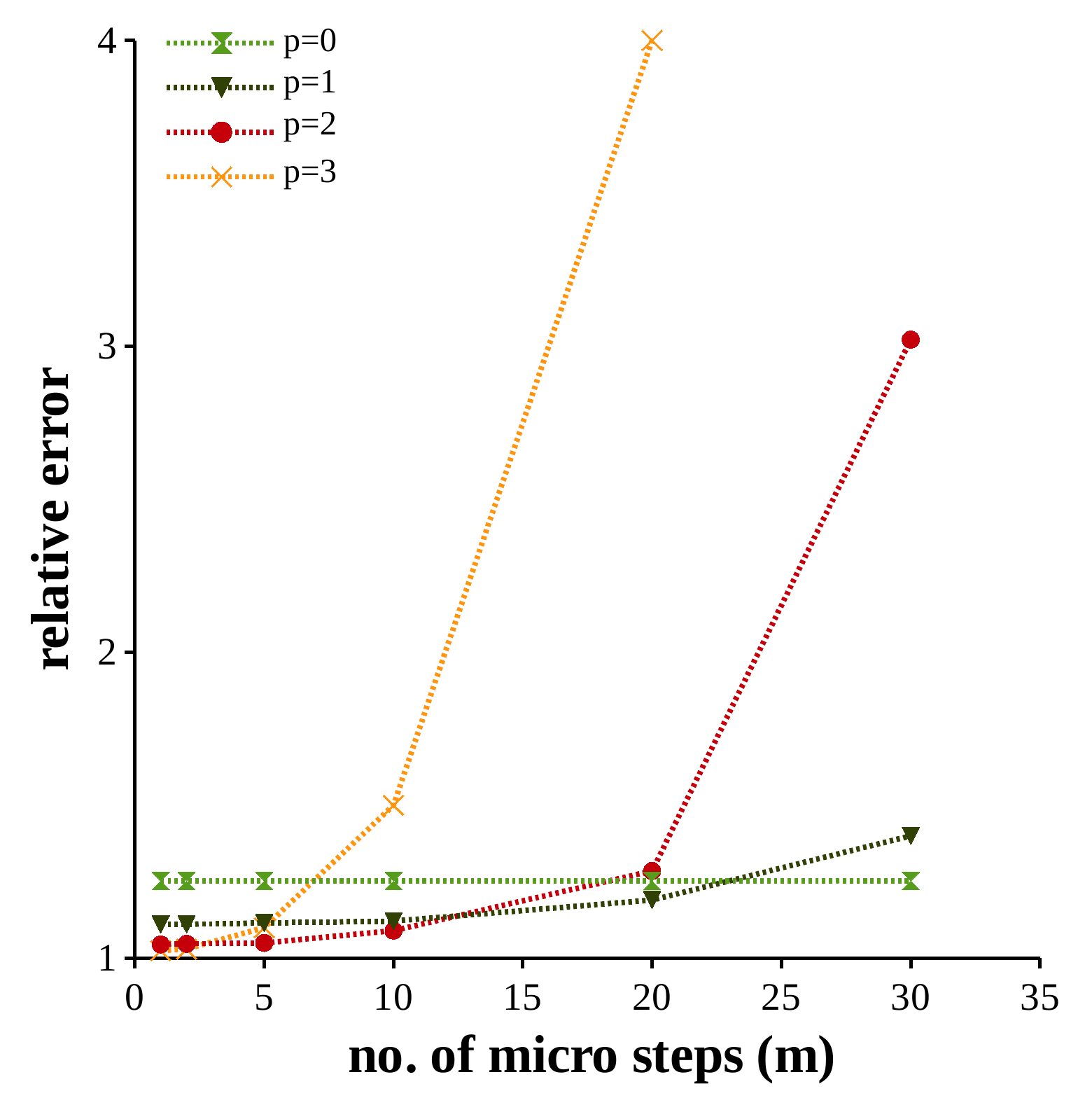}
  \caption{$1D$ test problem: Relative error in $P_g$ over $m$ at $t=18000$ sec for the semi-implicit fast-first MRT method using polynomials of order $p=0,1,2,3$ to approximate $\myvec{X}_G$ in the extrapolation steps. }
  \label{fig:ALG1-Pg-over-m-P0,P1,P2,P3ex}
\end{figure}

\subsubsection*{Speed up: }
Consider integration by the iteratively coupled and the semi-implicit MRT schemes between the synchronization level $n-1$ to $n$.  
There are $m$ uniform micro steps of size $h$ and one macro step of size $H=mh$ between $n-1$ and $n$.  
The iteratively coupled scheme solves both the flow and the geomechanical systems on the micro grid, 
while the semi-implicit MRT solves the flow system on the micro grid and the geomechanical system on the macro grid.  

Let $W_g$ be the time required to solve the geomechanical system, and $w_f$ be the time required for executing one Newton step in the flow system.
If the iterative scheme requires $n_{it}$ Newton steps to converge, then time required to solve the flow system per step is $W_f = n_{it} w_f$.
Further, if $n_{fp}$ fixed-point iteration steps are required to get the solution, then
the total time required for integration from $n-1$ to $n$ using the iteratively coupled scheme is, 

\begin{align*}
  \left(\text{CPU-time}\right)_{it} = m\ n_{fp} \left(W_f + W_g \right) \ .
\end{align*}

Similarly, if the semi-implicit MRT scheme takes $n_{mrt1}$ Newton steps to converge, 
then the time required to solve the flow system per micro step is,

\begin{align*}
 W_{f,mrt1} = n_{mrt1} w_f = n_{s1} W_f \ ,
\end{align*}

where, $n_{s1} = \dfrac{n_{mrt1}}{n_{it}}$. \newline
The total time required for integration from $n-1$ to $n$ using the semi-implicit MRT scheme is, thus, 

\begin{align*}
 \left(\text{CPU-time}\right)_{mrt1} = m \ n_{s1}\ W_f + W_g
\end{align*}

Therefore, the speed up is given as,

\begin{align} \label{eqn:speedup-ALG1}
 \text{speed-up} &= \frac{\left(\text{CPU-time}\right)_{it}}{\left(\text{CPU-time}\right)_{mrt1}} 
		  = \frac{m\ n_{fp} \left(W_f + W_g \right)}{m\ n_{s1}\ W_f + W_g} \notag\\
		 &= \frac{ n_{fp}\ m }{ n_{s1} \left(1-C\right) m + C }
\end{align}

where, $C = \dfrac{W_g}{W_f+W_g}$. \newline
In Eqn. \eqref{eqn:speedup-ALG1}, $0<C<1$, $n_{fp}\geq 1$, and $n_{s1} = 1$ for a stable system and $n_{s1}>1$ for an unstable system.

\paragraph{}
We can identify the following two special cases:

\begin{enumerate}
 \item For a stable system, for the case of $m=1$ (i.e. decoupled sequetial scheme) we get the minimum value of speed-up, 
    \begin{align}\label{eqn:ALG1-speedup-min}
    (\text{speed-up})_{min}= n_{fp} \ . 
    \end{align}
 \item For infinitely large $m$ (i.e. the limiting case of $m\rightarrow \infty$) we get the maximum value of speed-up, 
    \begin{align}\label{eqn:ALG1-speedup-max}
    (\text{speed-up})_{max}= \dfrac{n_{fp}}{n_{s1} \ \left(1-C\right)} \ .
    \end{align}
\end{enumerate}

\paragraph{}
% This behaviour is reflected very well in Fig. \ref{fig:ALG1-speedup-over-m-P0,P1,P2,P3ex} showing the numerically obtained speed-up curves for the extrapolation based semi-implicit MRT schemes using polynomials of order $p=0,1,2$.
The speed-up curves for the extrapolation based semi-implicit MRT schemes using polynomials of order $p=0,1,2$ are plotted in Fig. \ref{fig:ALG1-speedup-over-m-P0,P1,P2,P3ex}.
For comparison, the speed-up curve from Eqn. \eqref{eqn:speedup-ALG1} is also plotted in Fig. \ref{fig:ALG1-speedup-over-m-P0,P1,P2,P3ex} using $W_g = 0.9$ sec, $W_f=1.21$ sec, $n_{fp}=2$ and $n_{s1}=1$.
We can see that our numerical results for the speed-up reflect very well the behaviour expected from Eqn. \eqref{eqn:speedup-ALG1}.
The speed up curves for polynomial order $p=0,1,2$ coincide for those values of $m$ where the respective schemes are stable. 
The maximum speed-up for polynomial order $p=2$ is slightly lower due to instabilities at higher $m$, as predicted.

\paragraph{}
It can be further inferred that the higher the value of $C$, the higher is the speed-up.
This implies that for problems where the solution of the latent system is more time-consuming, for example in the non-linear case, a higher speed-up can be expected from the extrapolation based semi-implicit MRT schemes.  

\begin{figure}
 \centering
  \includegraphics[scale=0.5]{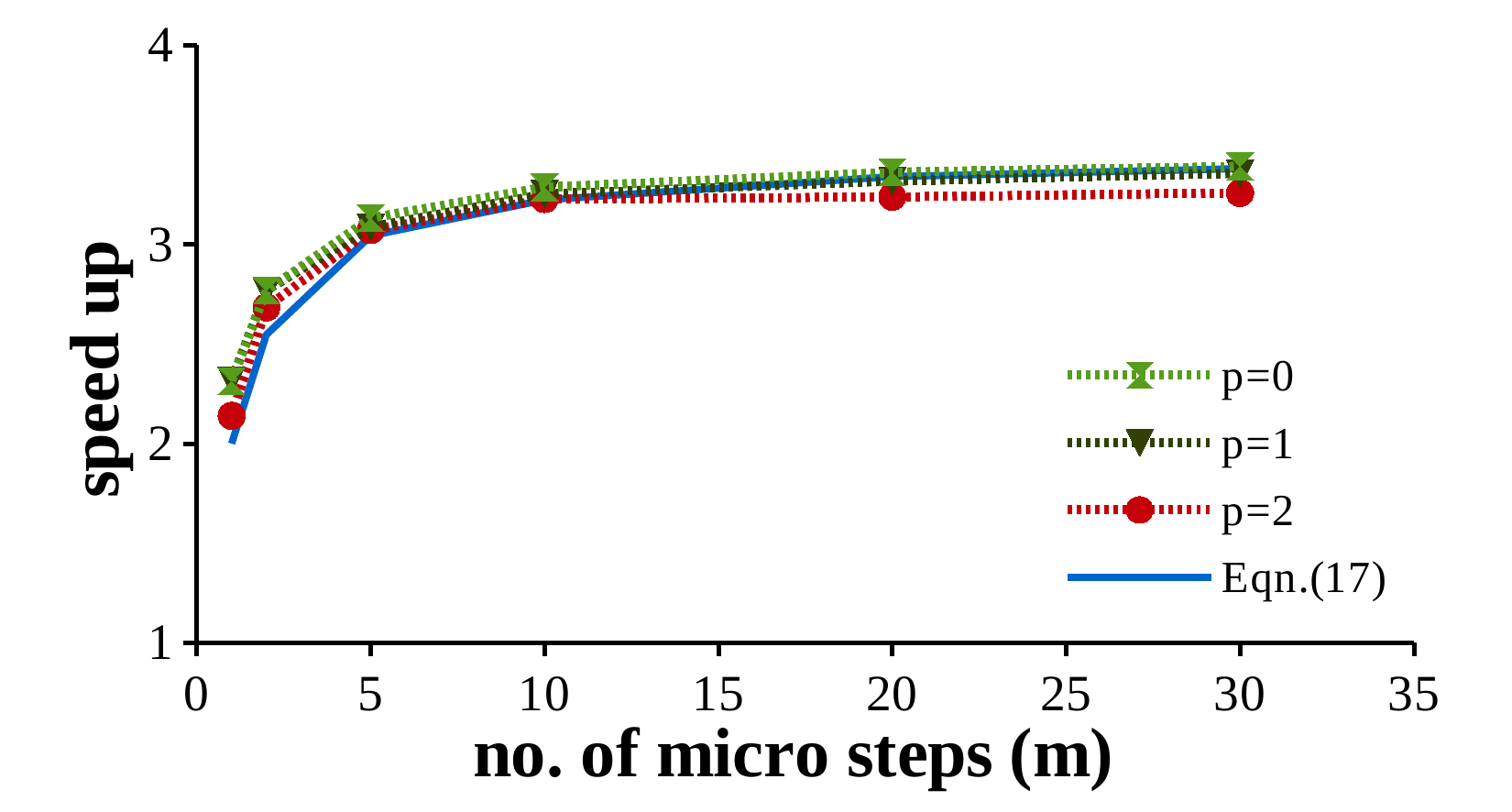}
  \caption{$1D$ test problem: Speed up in $P_g$ over $m$ for the semi-implicit fast-first MRT method using polynomials of order $p=0,1,2$ to approximate $\myvec{X}_G$ in the extrapolation steps. }
  \label{fig:ALG1-speedup-over-m-P0,P1,P2,P3ex}
\end{figure}

\subsubsection*{Stability of the scheme with extrapolation polynomial of order $p=2$: } 
\begin{figure}
 \centering
 \subfloat[$P_g$ profile over $z$ at $t=3600$ sec.]{
	      \includegraphics[scale=0.375]{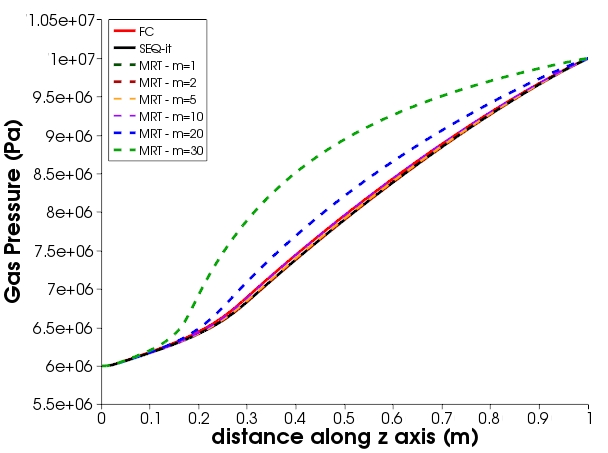} \label{fig:ALG1-Pg-over-z-wrt-m}
	    }
  \vfill
  \subfloat[$u_z$ profile over $z$ at $t=3600$ sec.]{
	      \includegraphics[scale=0.375]{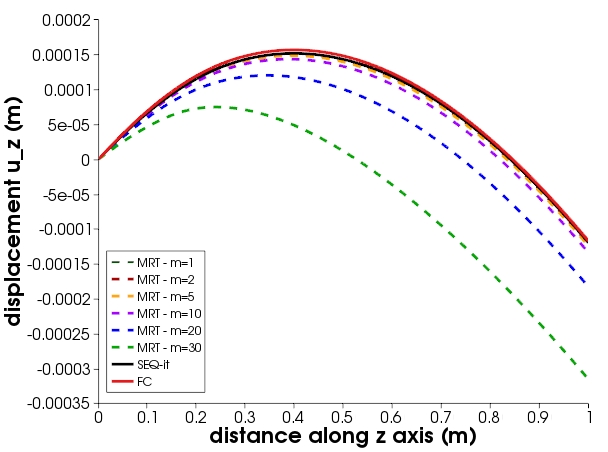} \label{fig:ALG1-uz-over-z-wrt-m}
	    }
  \caption{$P_g$ and $u_z$ profiles for the $1D$ test problem for fully-coupled (FC) scheme, iteratively-coupled (SEQ-it) scheme, and $m$-step semi-implicit MRT schemes using polynomial of order $p=2$ for extrapolation.}
  \label{fig:ALG1-PgUz-over-z-wrt-m}
\end{figure}

For problems where a high accuracy is desired, we usually use an extrapolation polynomial of order $p=2$, but this results in instabilities.
In Fig. \ref{fig:ALG1-Pg-over-z-wrt-m} and Fig. \ref{fig:ALG1-uz-over-z-wrt-m} showing the $P_g$ and the $u_z$ profiles, we can see that the solutions deviate from the reference solution very much for $m >10$ for this example.
This implies that we cannot choose any arbitrary value of $m$. 
The value of $m$ should be small enough so that the extrapolation error does not dominate, but it should be large enough so that we can take advantage of the speed up optimally, 
thus making the choice of $m$ an important consideration.

\subsubsection*{Choice of $m$: }
Consider the non-dimensional numbers $C_v$ and $C_r$.
$C_v$ is the consolidation coefficient which comes from Terzaghi's consolidation theory, and is given by,

$$ C_v = \frac{\kappa \left( \frac{k_{r,w}}{\mu_w} + \frac{k_{r,g}}{\mu_g} \right)}{\alpha_{biot}^2 K_m + S}\ .$$

$C_r$ is the reaction coefficient indicative of the damping of the normal consolidation due to dissociation kinetics, and is given by,

$$ C_r = \frac{\left(\frac{M_g}{\rho_g} + N_h \frac{M_w}{\rho_w} - \frac{M_h}{\rho_h}\right) k^{r} A_{r,s} }{\alpha_{biot}^2 K_m + S}\ ,$$

where, $\phi_e=\phi \left( 1-S_h\right)$ is the effective porosity, $S$ is the bulk storativity of the system given by, \newline
$S=\phi_e \left( K_w S_{w,e} + K_g S_{g,e} \right) + \left( \alpha-\phi_e\right)K_{sh} $, \
$K_{\gamma}$ is the compressibility of material $\gamma$, and $K_m=\dfrac{K_{sh}}{1-\alpha_{biot}}$ is the compressibility of the bulk porous material.

The derivation of $C_r$ can be found in \cite{Gupta2015} for a simplified 1D consolidation problem with hydrate dissociating in the sample due to depressurization.

\paragraph{}
The rate of consolidation is directly proportional to $C_v$ and inversely proportional to $C_r$. 
Thus, the $C_v/C_r$ ratio can be seen as the relative activity of the latent component.

\begin{figure}
 \centering
  \includegraphics[scale=0.5]{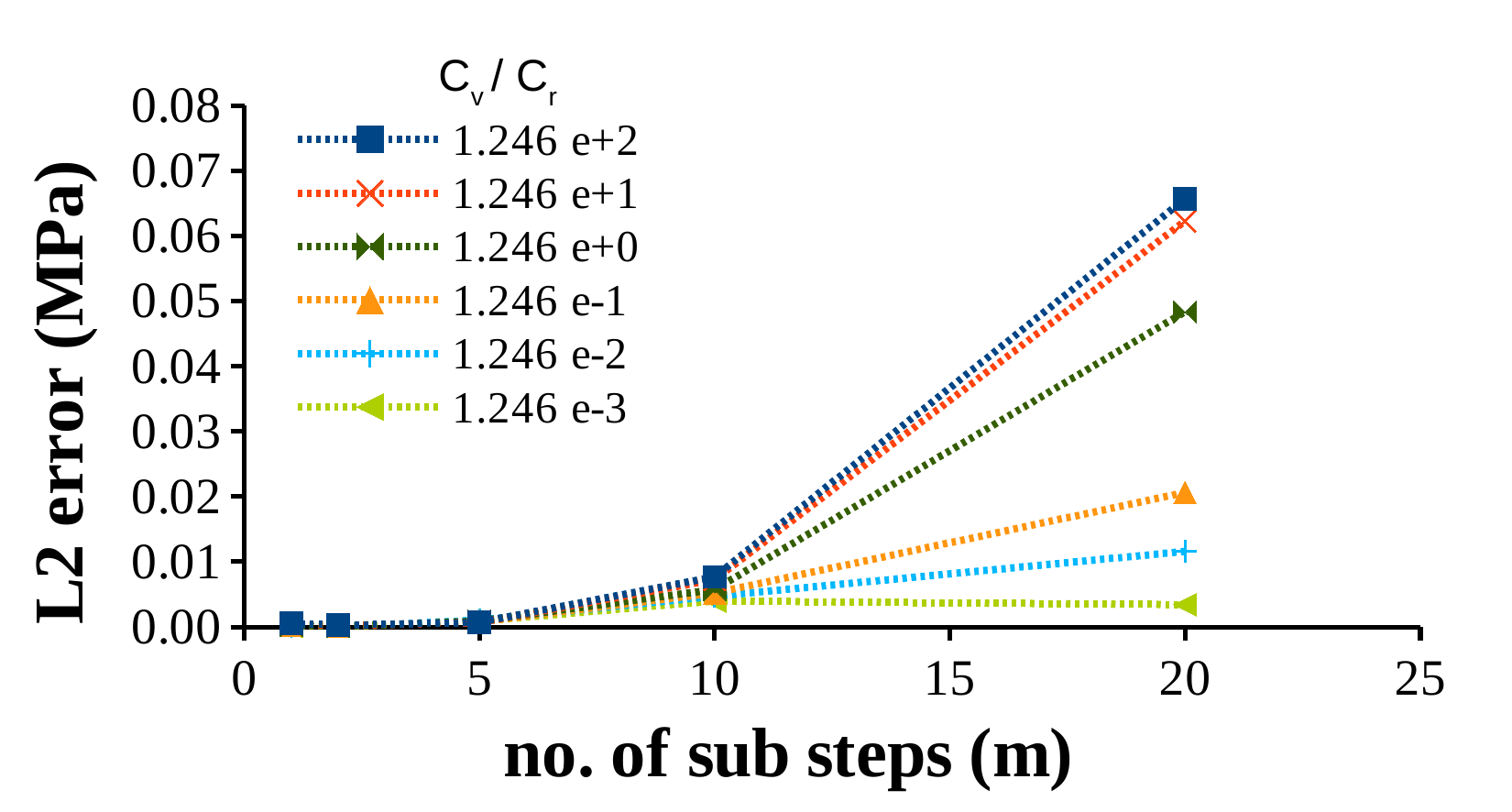}
  \caption{$1D$ test problem: Error in $P_g$ plotted over $m$ for different values of the $C_v/C_r$ ratios for $m$-step semi-implicit MRT schemes using polynomial of order $p=2$ for extrapolation.}
  \label{fig:ALG1-L2errPg-CvbyCr}
\end{figure}

\paragraph{}
In Fig. \ref{fig:ALG1-L2errPg-CvbyCr}, we show the errors in $P_g$ plotted over $m$ for a broad range of $C_v/C_r$ ratios.
We observe that the higher the relative activity of the latent component, the more dominating is the extrapolation error.
For the case of a non-dominating extrapolation error, the values of $m$ can be chosen very large. 
However, in the absence of a priori estimate of the relative activity, the value of $m$ must be kept small in order to keep the relative error close to $1$.
In our hydrate reservoir simulator we choose $m\leq 5$.  

% \begin{figure}
%  \centering
%   \includegraphics[scale=0.5]{ALG1-SpeedUp_differentCvCrRatio.pdf}
%   \caption{Please write your figure caption here}
%   \label{fig:ALG1-speedup-CvbyCr}
% \end{figure}

\subsubsection{Results - Compound-fast MRT method}
\subsubsection*{Relative error: }
In Fig. \ref{fig:ALG2-L2errPg-vs-m} and Fig. \ref{fig:ALG2-L2errU-vs-m}, we show the relative errors in $P_g$ and $u_z$ plotted over $m$, respectively, for the compound fast MRT method.
For comparison, we have also plotted the errors in $P_g$ and $u_z$ for the semi-implicit MRT method using extrapolation polynomials of order $p=2$.
We can see that the relative errors for this scheme remain close to $1$ and do not depend on $m$ significantly.

\paragraph{}
In this example, the active ODE becomes unsolvable in the predictor step at $m\geq 60$ for the minimum error reduction of $10^{-3}$ that is prescribed for the Newton solver.

\subsubsection*{Speed-up: }
\begin{figure}
 \centering
 \subfloat[Relative error in $P_g$ over $m$]{
	    \includegraphics[scale=0.5]{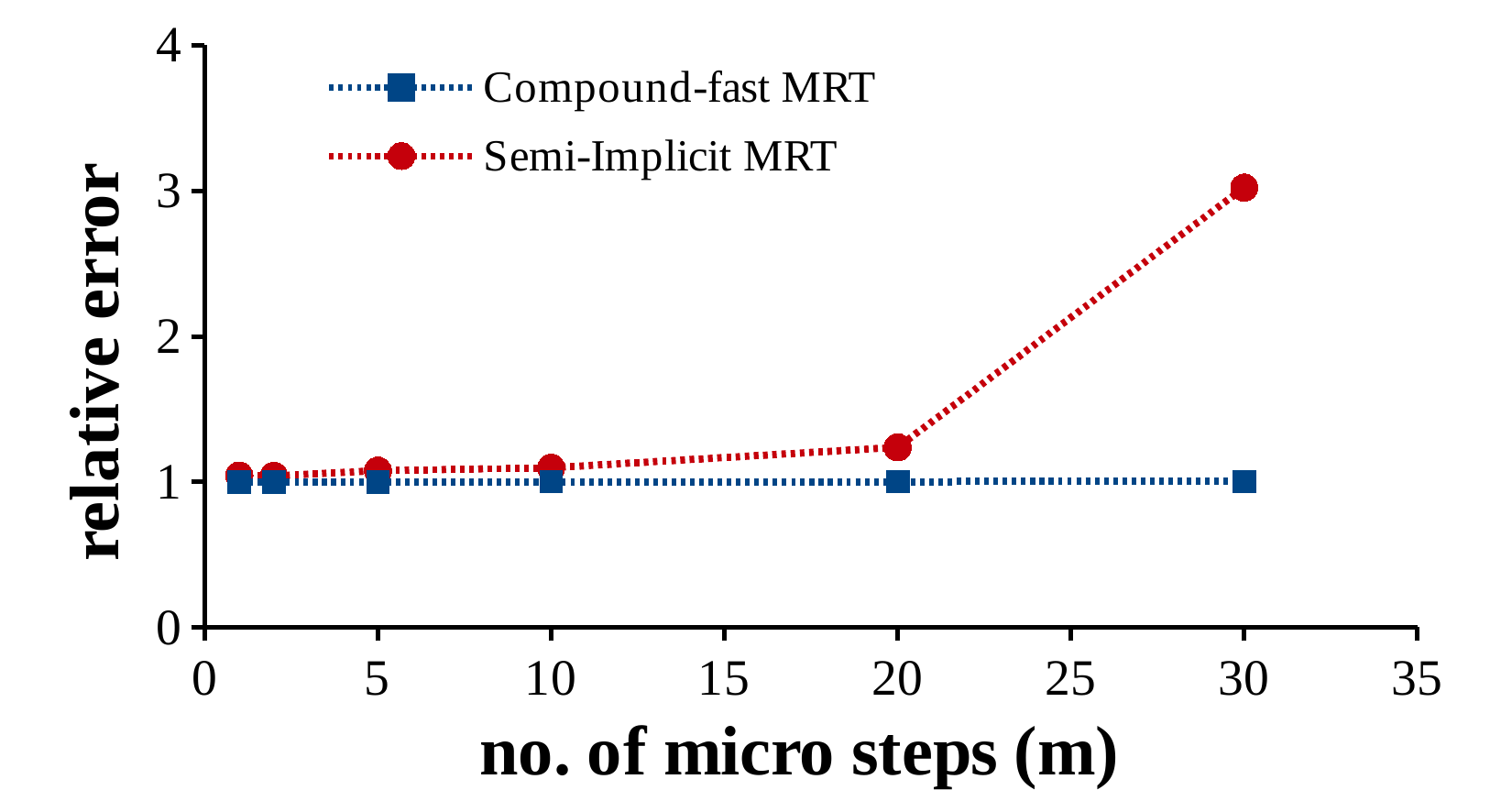} \label{fig:ALG2-L2errPg-vs-m}
	  }
  \vfill
  \subfloat[Relative error in $u_z$ over $m$]{
	    \includegraphics[scale=0.5]{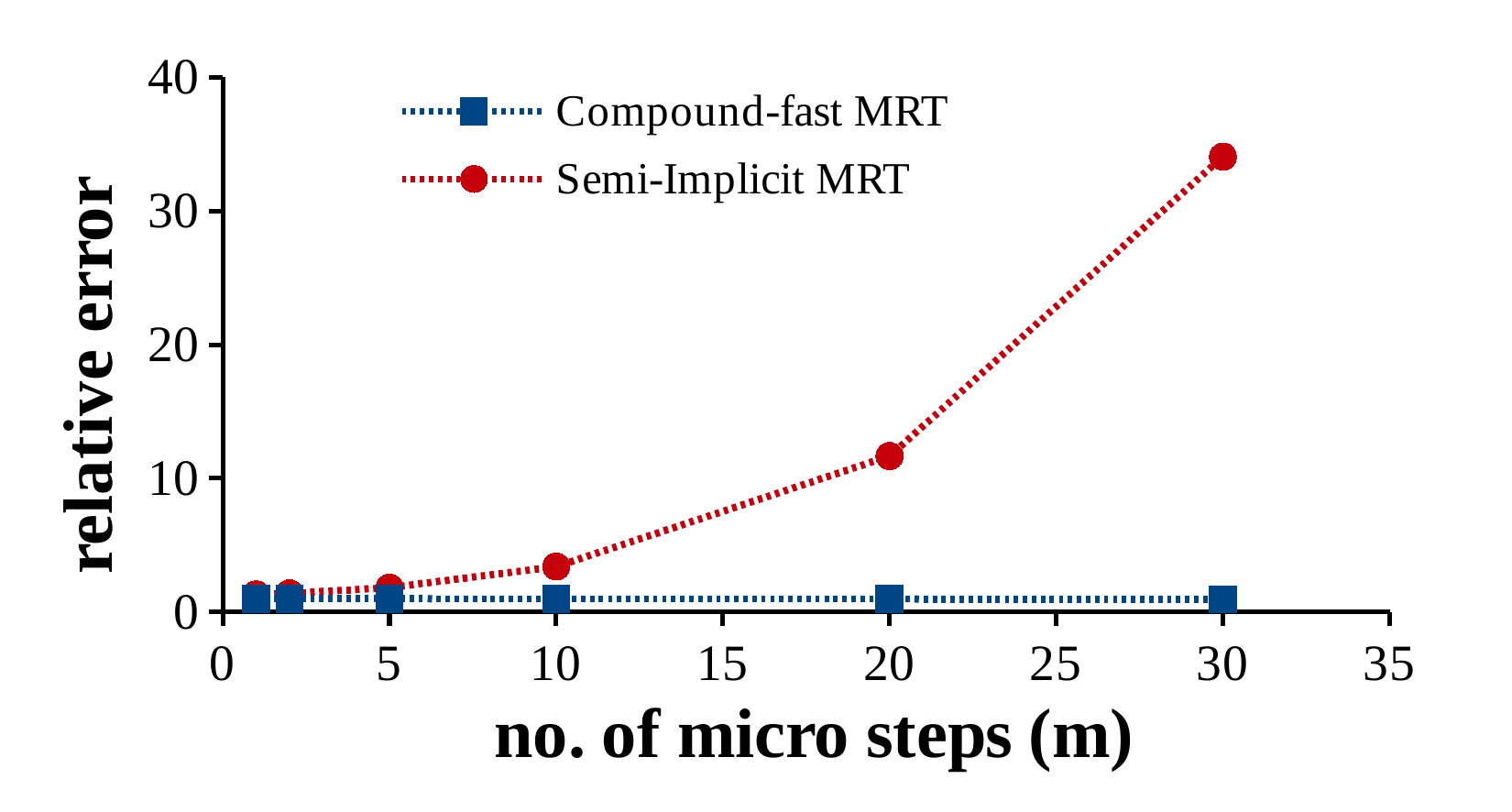} \label{fig:ALG2-L2errU-vs-m}
	  }
  \vfill
   \subfloat[Speed-up over $m$]{
	    \includegraphics[scale=0.5]{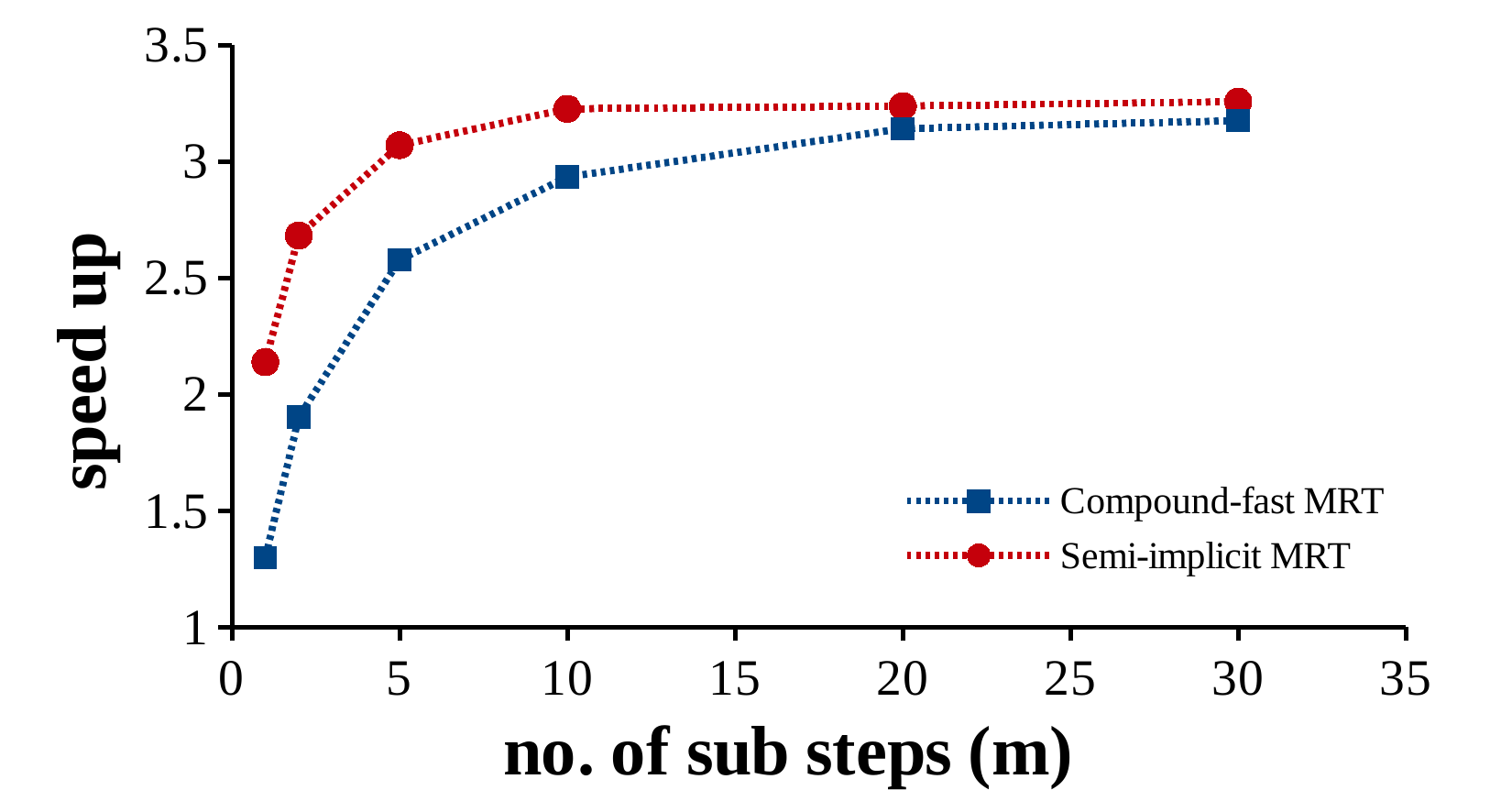} \label{fig:ALG2-speedup-vs-m}
	  }
  \caption{$1D$ test problem: Relative error and speed-up curves for the compound-fast MRT method and the semi-implicit MRT method using polynomial of order $p=2$ for extrapolation.}
\end{figure}

Fig. \ref{fig:ALG2-speedup-vs-m} shows the speed up obtained from this scheme as compared to the speed up from the semi-implicit scheme.
We can see that this scheme gives a lower speed up, especially for smaller values of $m$.

\paragraph{}
The speed-up curve can be derived for the compound-fast method using similar arguments as in the case of the semi-implicit MRT method. 
The compound-fast method solves the flow system once on the macro-grid (during the predictor step) and once on the micro grid ($m$ micro-steps), 
and the geomechanical system twice on the macro grid (once during the predictor step, and once during the corrector-macro step).
If the compound-fast scheme takes $n_{mrt2,p}$ and $n_{mrt2}$ Newton steps to converge for the predictor-step and for each micro-step, respectively, 
then the time required to integrate the flow system per micro step is,

\begin{align*}
 W_{f,mrt2} &= n_{mrt2}\ w_f = n_{s2} \ W_f \ ,
\end{align*}

and the time required to solve the flow system for the predictor step is,

\begin{align*}
 W_{f,p} = n_{mrt2,p}\ w_f = n_{s2,p} \ W_f \ ,
\end{align*}

where, $n_{s2,p} = \dfrac{n_{mrt2,p}}{n_{it}}$, and $n_{s2} = \dfrac{n_{mrt2}}{n_{it}}$.
Furthermore, since the scheme is stable, $n_{s2} = 1$. \newline
The total time required for integration from $n-1$ to $n$ using the compound-fast MRT scheme is, thus, 

\begin{align*}
 \left(\text{CPU-time}\right)_{mrt2} = W_{g,p} + \left( n_{s2,p} + m \right) W_f + W_{g}
\end{align*}
where, $W_{g,p}$ is the time required to solve the geomechanical system in the predictor step. 

\paragraph{}
Therefore, the speed-up is given as,

\begin{align} \label{eqn:speedup-ALG2}
 \text{speed-up} &= \frac{\left(\text{CPU-time}\right)_{it}}{\left(\text{CPU-time}\right)_{mrt2}} \notag\\
		 &= \frac{m\ n_{fp} \left(W_f + W_g \right)}{W_{g,p} + \left( n_{s2,p} + m  \right) W_f + W_{g}} \notag\\
		 &= \frac{ n_{fp}\ m }{ \left(1-C\right) m + C + \Delta_p}
\end{align}

where, $\Delta_p = n_{s2,p} \left( 1-C\right) + C_p$ , s.t., $C_p = \dfrac{W_{g,p}}{W_f + W_g}$ . \newline
In Eqn. \eqref{eqn:speedup-ALG2}, $\Delta_p >0$ because $n_{s2,p} > 0$, $C_p > 0$ and $0<C<1$.

\paragraph{}
Comparing equations \eqref{eqn:speedup-ALG1} and \eqref{eqn:speedup-ALG2}, we can conclude that for any given $m\geq1$, 
the speed-up obtained from the compound-fast MRT method is always smaller than that obtained from the corresponding semi-implicit MRT method.
However, for the limiting case of $m\rightarrow \infty$, the speed-up from the compound-fast method approaches the speed-up from a stable semi-implicit method. 

\subsubsection{Discussion}
\paragraph{}
Based on the above results, we now compare the performance of the semi-implicit and the compounf-fast MRT methods, and
discuss the factors which affect the choice of a particular MRT method for a given problem.

\subsubsection*{Advantages of the Semi-implicit MRT method: }

\begin{itemize}
 \item This scheme gives a higher speed up for a given $m$, especially if the natural time scale of the latent component is comparable to that of the active component.    
%  \item This scheme is more attractive when the solution of the latent system is more time-consuming, for example in case of a $3D$ problem where the matrix assembly times are larger, 
% 	s.t., $\Delta_p$ is large.
% \item This scheme can be more attractive when the work ratio, $C$, is higher (See Eqn. \eqref{eqn:ALG1-speedup-max}).
\end{itemize}

\subsubsection*{Disadvantages of the Semi-implicit MRT method: } 

\begin{itemize}
  \item The scheme with lower order polynomilal for extrapolation shows high stability but low accuracy. 
	The use of higher order polynomial extrapolation improves the accuracy significantly, but introduces errors which grow exponentially with $m$, thus, making the scheme unstable at large $m$ values.
  \item With the use of higher order polynomials, the choice of $m$ becomes strongly dependent on the actual activity of the latent components.
%   \item Since the semi-implicit MRT method is essentially a decoupled sequential solution strategy, 
% 	it places more severe restrictions on the requirement of weak coupling as compared to the compound-fast MRT method \cite{Verhoeven2007}.
  \item For problems where relative activities of the components are expected to fluctuate over time, one must use MRT scheme which is relatively insensitive to the component activities.
	Therefore, for such problems the semi-implicit scheme can only be used with lower order polynomial extrapolation, thus compromising the accuracy of the solution. 
\end{itemize}

\subsubsection*{Advantages of the Compound-fast MRT method: }

\begin{itemize}
 \item This scheme is stable for arbitrarily large values of $m$, provided that the active system is stable and solvable upto the prescribed error reduction in the predictor-step \cite{Verhoeven2007}.
 \item The errors in the solution are comparable to that of the iteratively coupled solution scheme and do not grow with increasing $m$.
  The time-mesh for the latent component can be made quite coarse irrespective of the actual activity of the latent component.
 \item It is easier to handle problems where the activity of latent system fluctuates in time
  because the stability of the scheme is insensitive to the refinement of the macro-grid. 
%  \item It is also possible to adaptively control the macro time-steps using the number of Newton steps made in the predictor step, as described in Algorithm 3. 
\end{itemize}

% This scheme is stable for arbitrarily large values of $m$, provided that the active system is stable and solvable upto the prescribed error reduction in the predictor-step. 
% The errors in solution are comparable to that of the iteratively coupled solution scheme and do not grow with increasing $m$.
% The value of $m$ can thus be chosen to reflect the actual activity of the latent components, which is certainly attractive for problems with a very large difference in time scales.
% It is easier to adaptively control the time integration using this scheme because the errors in solution do not depend on $m$ and the choice of $m$ itself does not affect the stability of the scheme.
% It is sufficient to ensure the solvability of the of the active ODE in the predictor-step. 

\subsubsection*{Disadvantages of the Compound-fast MRT method: } 

\begin{itemize}
  \item It gives lower speed-up as compared to the semi-implicit scheme, especially for lower values of $m$.
  The difference in the speed-ups between the two MRT schemes becomes more pronounced for more complex examples. 
  It mostly depends on the work ratio $C$ and the work ratio in the predictor step $C_p$ .
  The general trend is, the closer to $1$ the ratio $C / C_p$ is, the larger is the difference in the speed-ups.
  (Compare Eqns. \eqref{eqn:speedup-ALG1} and \eqref{eqn:speedup-ALG2}.)
  \item
  If the natural time scale of the latent system is large, then the increase in computation time due to predictor steps can be compensated by using a very coarse macro grid. 
  If, however, the natural time scale of the latent system is small and/or the time-evolution of the latent component is desired at a finer resolution, then this method may prove computationally more intensive.
\end{itemize}

% It gives lower speed up as compared to the semi-implicit scheme, especially for lower values of $m$, which is expected because of the extra steps made in the predictor step.
% This can infact get very time-consuming if the latent system takes more time to solve. 
% If the natural time scale of the latent system is large, then this increase in computation time due to predictor steps can be compensated by the coarsened macro grid. 
% If, however, the natural time scale of the latent system is small, this method may prove computationally more intensive.

\subsection{Test 2}
\paragraph{}
We now show the performance of the MRT schemes for a relatively complex example in a $3D$ setting 
where, a sub-surface hydrate reservoir is destabilized through depressurization using a vertically placed low pressure gas well. 

\subsubsection{Problem setting}
\begin{figure}[h]
    \centering
      \includegraphics[scale=0.4]{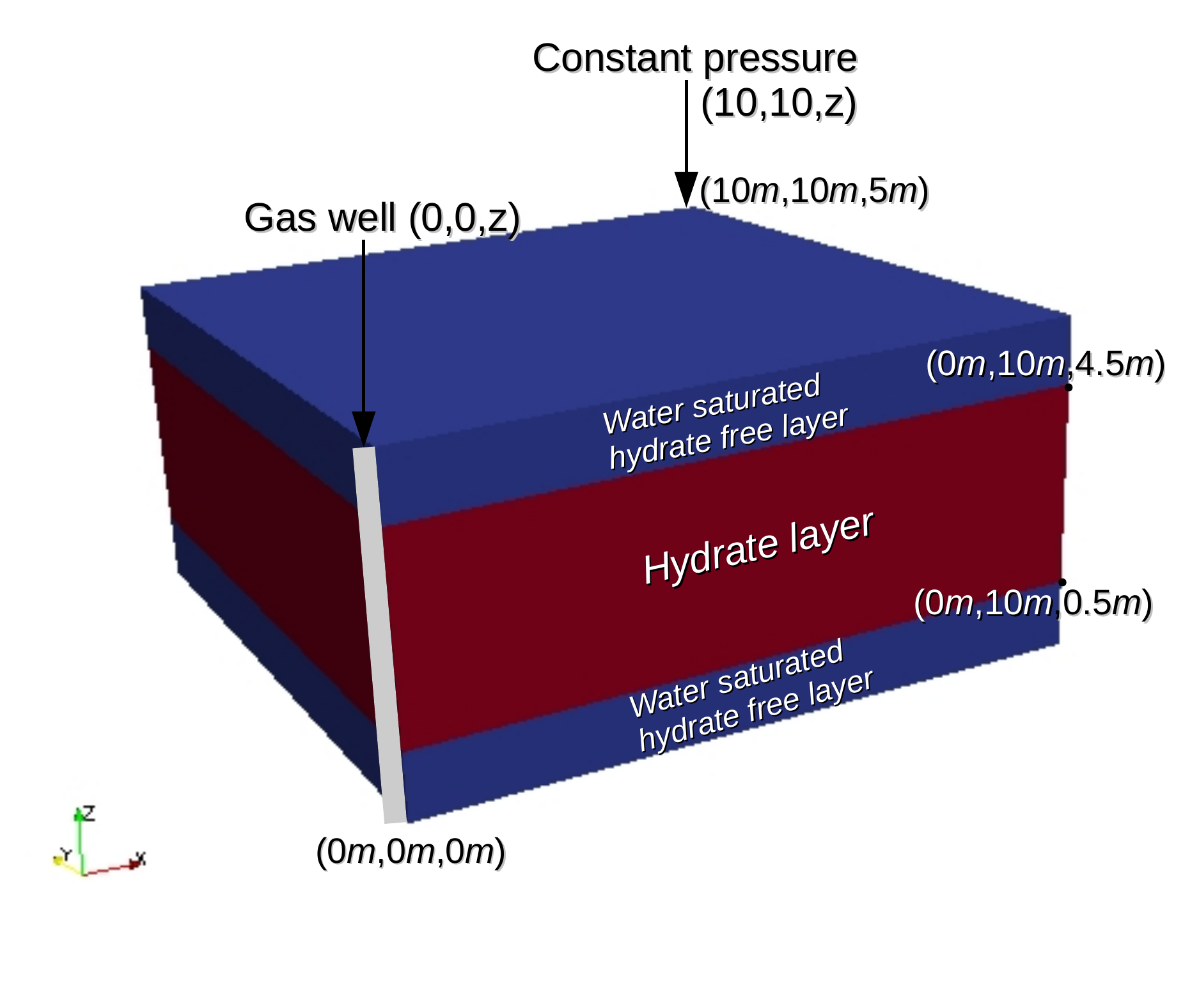}
      \captionof{figure}{Schemmatic of the $3D$ hydrate reservoir problem}
      \label{fig:3DtestSchemmatic}
  \end{figure}

\paragraph{}
We consider a scaled down $3D$ reservoir with dimensions $10m\times10m\times5m$, as shown in Fig. \ref{fig:3DtestSchemmatic}.
The hydrate is homogeneously distributed in a $4m$ thick layer lying between $0.5m\leq z \leq 4.5m$, and has a saturation of $40\%$ by volume. 
The reservoir is fully saturated with water and has an initial pressure of $10$ MPa. 
The reservoir is depressurized through a low pressure gas well located at $(0,0,z)$.
The pressure in the gas well is maintained at $P_{well}=4$ MPa.
A constant vertical load of $10$ MPa is acting on the top boundary of the reservoir (i.e. at $z=10$ m).
The initial and the boundary conditions are listed in Table \ref{table:3Dtest_ICValues} and Table \ref{table:3Dtest_BCValues} respectively.
% The material properties and other model parameters are listed in Table \ref{table:3Dtest_properties}.

  \begin{table}
    \centering
    \captionof{table}{Initial conditions for the $3D$ hydrate reservoir problem }\label{table:3Dtest_ICValues}
   \begin{tabular}{L{4cm} L{2.5cm}}	\hline \\ [-0.7em]
    \textbf{Hydrate layer} 				\\
						& $P_{eff,i} = 10$ MPa 		\\
						& $S_{h,i} = 0.4$ 		\\
    at $t=0$, for				& $S_{g,i} = 0$			\\ 
    $0\leq x,y \leq 10$ , $0.5\leq z \leq 4.5$	& $T_{i} = 10 \ ^0$C		\\
						& $K_{i} = 0.0198$ mD		\\
						& $\phi_{eff,i} = 0.18$		\\ [0.4em] \hline \\ [-0.7em]
    \textbf{Hydrate-free layers} \\
						& $P_{eff,i} = 10$ MPa 		\\
    at $t=0$, for				& $S_{h,i} = 0$ 		\\
    $0\leq x,y \leq 10$ , $z < 0.5$		& $S_{g,i} = 0$			\\ 
    and						& $T_{i} = 10 \ ^0$C		\\
    $0\leq x,y \leq 10$ , $z > 4.5$					& $K_{i} = 0.1$ mD		\\
						& $\phi_{eff,i} = 0.3$		\\ [0.4em] \hline    
   \end{tabular}
  \end{table}
  
  \begin{table}
   \centering
   \captionof{table}{Boundary conditions for the $3D$ hydrate reservoir problem }\label{table:3Dtest_BCValues}
   \begin{tabular}{L{4.5cm} L{2.25cm}}
    \multicolumn{2}{c}{FLOW model} 						\\ \hline
    Gas well at 			& $P_g = 4$ MPa 			\\
    $x=0$, $y=0$, $0\leq z \leq 5$	& $S_w = 0$ 				\\
				 	&$\nabla\cdot T = 0$			\\ \hline
    Pressure constraint at		& $P_{eff} = P_{eff,i}$ 		\\
    $x=10$, $y=10$, $0\leq z \leq 5$	& $S_w = S_{w,i}$ 			\\
					& $T = T_i$				\\ \hline
    No-flow and adiabatic conditions 	& $\mathbf{v}_g \cdot  \hat n = 0$ 	\\
    on remaining boundaries, i.e.,	& $\mathbf{v}_w \cdot  \hat n = 0$ 	\\
					& $\nabla\cdot T = 0$			\\ \hline      
    \\
    \multicolumn{2}{c}{GEOMECHANICAL model} 								\\ \hline
    \multicolumn{2}{l}{Top boundary} 									\\
    $0\leq x\leq 10$ , $0\leq y\leq 10$ , $z=5$	& $\sigma_{zz}=10$ MPa , $\sigma_{xy}=\sigma_{yx}=0$	\\ 
    \multicolumn{2}{l}{Bottom boundary} 								\\
    $0\leq x\leq 10$ , $0\leq y\leq 10$ , $z=0$ & $u_z=0$ , $\sigma_{xy}=\sigma_{yx}=0$			\\
    \multicolumn{2}{l}{Remaining boundaries}			 					\\ 
						& $u_x=u_y=0$ , $\sigma_{zz}=0$ 			\\ \hline
   \end{tabular}

  \end{table}

\subsubsection{Numerical simulation}

\paragraph{}
We discretize the space domain into $30\times30\times15$ cells
and simulate this problem using the compound-fast MRT method, and the semi-implicit MRT method with $p0$ polynomial extrapolation.
We assume uniform, non-adaptive micro and macro time-grids, and
chose the following values of the multirate factor,

$$
m = [1,2,5,10,15,30] \ .
$$
The simulation is run until $t_{end}=30000$ s.
The micro grid size is chosen as $h=200$ s, and the macro grid size as $H = m\ h$.

\subsubsection{Results}

\begin{figure*}
 \centering
    \subfloat[Hydrate saturation]{%
    \includegraphics[scale=0.23]{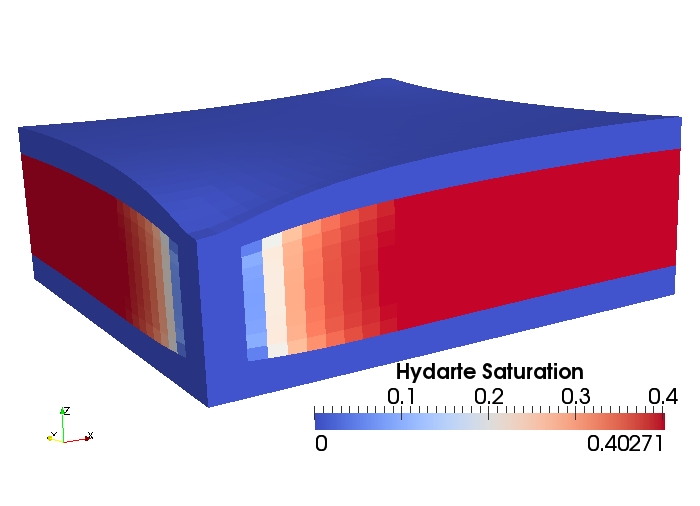}
    \label{fig:3Dtest-ShProfile}}
    \hspace{-0.5cm}
    \subfloat[Gas saturation]{%
    \includegraphics[scale=0.23]{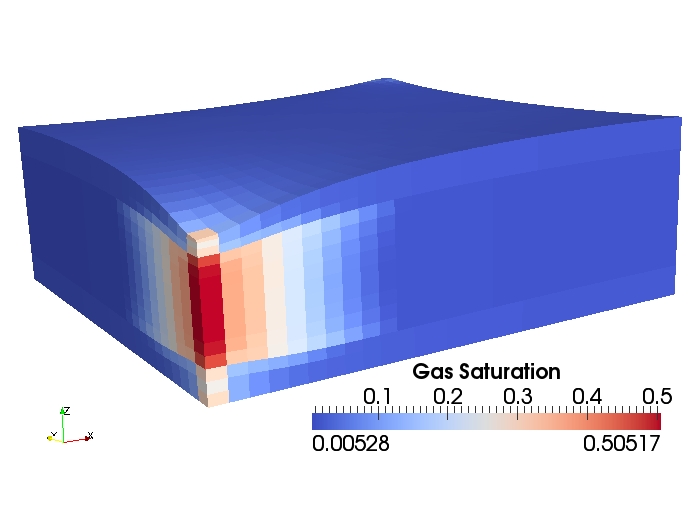}
    \label{fig:3Dtest-SgProfile}}
    \hspace{-0.5cm}
    \subfloat[Deviatoric stress]{%
    \includegraphics[scale=0.23]{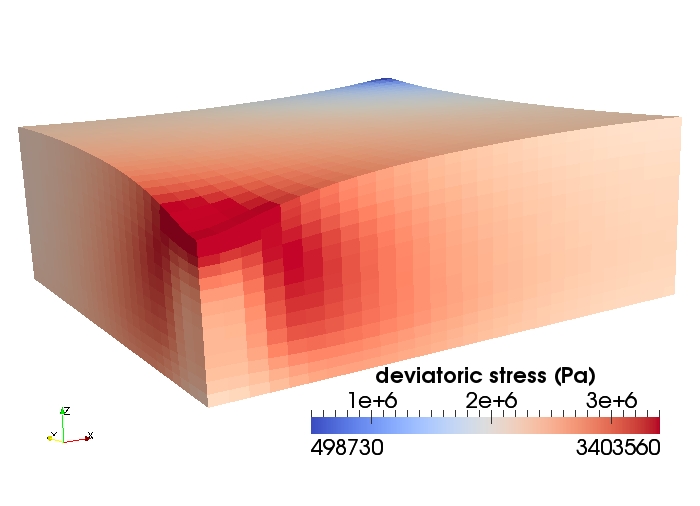}
    \label{fig:3Dtest-DevStressProfile}}
    \caption{$3D$ hydrate reservoir problem: Selected profiles at $time=t_{end}$.
    \\ The domain is warped with respect to displacement to show the ground subsidence around the well clearly. The warping of the domain is achieved through post-processing using PARAView.}
    \label{fig:3Dtest-results}
\end{figure*}

\paragraph{}
In this example, we simulate the melting of hydrate, methane gas generation, and the resulting ground subsidence and stress build-up in the vicinity of the well,
which are shown as screenshots at $t_{end}$ in Fig. \ref{fig:3Dtest-results}.

\begin{figure}[h]
 \centering
  \includegraphics[scale=0.5]{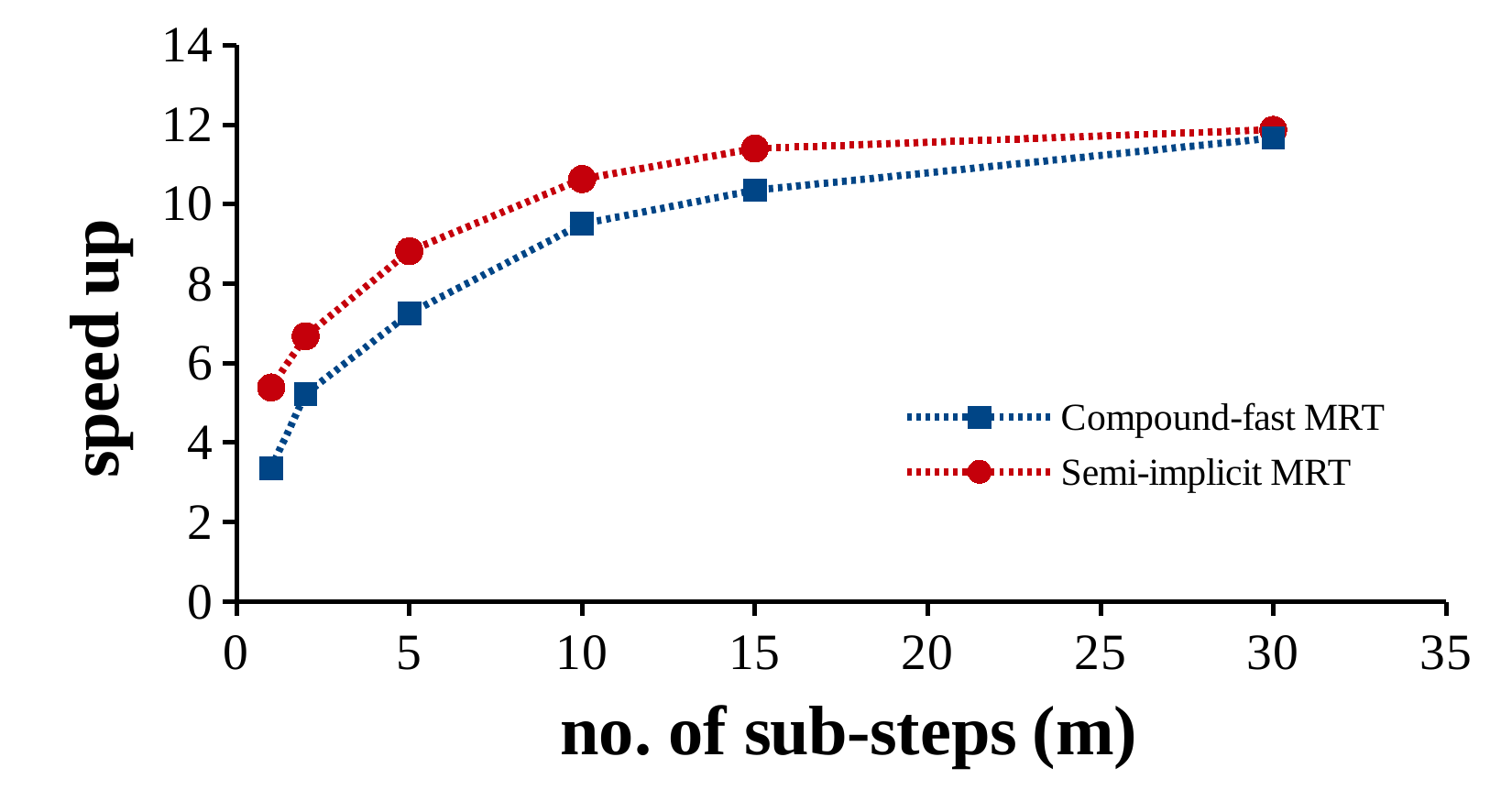}
  \caption[$3D$ hydrate reservoir problem: Speed up over $m$ for the compound-fast and the semi-implicit MRT methods]
  {$3D$ hydrate reservoir problem \\Speed up over $m$ for the compound-fast MRT method, and the semi-implicit MRT method using $p=0$ polynomial for extrapolation.}
  \label{fig:3Dtest-speedup}
\end{figure}

\paragraph{}
To evaluate the performance of the MRT methods in comparison to the iterative solution scheme, we obtain the $speed-up$ vs. $m$ curves for the semi-implicit and the compound-fast MRT methods.
which are shown in Fig. \ref{fig:3Dtest-speedup}. We obtain a maximum speed up of approximately $12$ for this example problem.
We can see that, by using the MRT methods for time-integration we can get a significant speed up for the $3D$ case.
% To give some perspective, for this problem the CPU-time required to obtain the solution at $t_{end}$ using compound-fast MRT method with $m=30$ was approximately $1$ day.
% In comparison, the CPU-time required by the iterative Gauss-Seidel scheme to obtain the solution with comparable accuracy was over $12$ days.

\subsubsection{Discussion}
\paragraph{}
% This large speed-up in the $3D$ case can be attributed to the fact that the solution of the latent system is more time-consuming in case of a $3D$ problem due to larger matrix assembly times.
This large speed-up in the $3D$ case can be attributed to the fact that the work ratio $C$ is higher due to the larger matrix assembly times in case of a $3D$ problem as compared to the $1D$ problem.
Consider the expression for maximum speed-up, Eqn. \eqref{eqn:ALG1-speedup-max}, where,
for the $1D$ test case, $C\approx 0.43$ while for the $3D$ test case, $C\approx 0.84$, i.e., almost twice as much as that for the $1D$ case.

\section{Conclusions}\label{sec:conclusions}
\paragraph{}
For multi time scale hydro-geomechanical subsurface flow problems, the multirate time stepping methods provide a significant speed up as compared to fully coupled or decoupled (iterative or sequential) schemes, 
especially for $3D$ problems, provided the model can be partitioned into sufficiently weakly coupled subsystems having distinctly different time scales. 
In our case, we deal with subsurface hydrate reservoirs where the mathematical model is naturally partitioned into the active flow system and latent geomechanical system.
% In this paper, besides outlining the MRT strategy that we use in our hydrate reservoir simulator, we also focus on the factors which affect the choice of a particular MRT method for a given problem.
The stability of the extrapolation-based semi-implicit method is sensitive to the activity of the latent component, 
while the compound-fast method is fairly independent of the activity of the latent component.
If the difference in the time scales between the active and latent components is comparable, then the semi-implicit MRT is more attractive as it gives a higher speed up.
It must, however, be kept in mind that this method is only conditionally stable for higher-order extrapolation and the extrapolation errors tend to accumulate with increasing $m$. 
It is therefore necessary to keep the choice of $m$ small. 
On the other hand, if the difference in the time scale between the active and latent components is large, then the compound-fast MRT method is more suitable as it is stable for arbitrarily large values of $m$, 
provided that the active system is solvable in the predictor step.
Another important consideration is whether the relative activities are expected to vary over time.
For problems where the activity of the latent component in particular fluctuates in time, it is important that the MRT method be insensitive to the the activity of the latent component, thus
making compound-fast MRT methods more attractive in such cases. 
% Also, if adaptive time step control is desired to increase the speed up,
% the compound-fast MRT method comes in handy as it is possible to heuristically control the macro time-steps by using the number of newton steps made in the predictor step, as described in Algorithm 3.

\paragraph{}
Many extensions of these MRT methods are possible in our hydrate reservoir simulator, for example, accounting for local variations in time scales of each component over the space domain, and,
stabilization of the semi-implicit MRT method, etc., which can make the MRT methods more attractive for solving large scale problems, especially in $3D$, more efficiently.

% \end{linenumbers}

% \section*{Acknowledgements}
% We gratefully acknowledge the support for the first author by the German Research Foundation (DFG) through project no. WO 671/11-1.

\vspace{0.5cm}





\end{document}